\documentclass[journal,12pt,onecolumn,draftclsnofoot,]{IEEEtran}

\pdfoutput=1
\IEEEoverridecommandlockouts                              % This command is only needed if 
\newtheorem{Remark}{Remark}
\newtheorem{Corollary}{Corollary}

\newenvironment{Proof}{\noindent{\em Proof:\/}}{\hfill $\Box$\par}
\newtheorem{Theorem}{Theorem}
\newtheorem{Lemma}{Lemma}

\newtheorem{Assumption}{Assumption}

\usepackage{graphicx}
\usepackage{amsmath,bm}
\usepackage{amssymb}
\usepackage{algorithm}
\usepackage{algpseudocode}
\usepackage{cite}
\usepackage{hyperref}
\usepackage[caption=false]{subfig}
\hypersetup{hidelinks}
\algdef{SE}[DOWHILE]{Do}{doWhile}{\algorithmicdo}[1]{\algorithmicwhile\ #1}%

\newcommand{\mathactivatecomma}{%
  \begingroup\lccode`~=`\,
  \lowercase{\endgroup\edef~}{\mathchar\the\mathcode`\,\penalty0 }}

\algnewcommand{\Initialize}[1]{%
  \State \textbf{Initialize: $i \in \mathcal{V}$}
  \Statex \hspace*{\algorithmicindent}\parbox[t]{.8\linewidth}{\raggedright #1}
}

\algnewcommand{\Iteration}[1]{%
  \State \textbf{Iteration $(t\geq 0)$: $i \in \mathcal{V}$}
  \Statex \hspace*{\algorithmicindent}\parbox[t]{.8\linewidth}{\raggedright #1}
}

\algnewcommand{\Output}[1]{%
  \State \textbf{Output: $i \in \mathcal{V}$}
  \Statex \hspace*{\algorithmicindent}\parbox[t]{.8\linewidth}{\raggedright #1}
}

\begin{document}

\bstctlcite{IEEEexample:BSTcontrol}

\title{A Gradient-Free Distributed Optimization Method for Convex Sum of Non-Convex Cost Functions}

\author{Yipeng Pang and Guoqiang Hu % <-this % stops a space
% \thanks{This research was supported by Singapore Ministry of Education Academic Research Fund Tier 1 RG180/17(2017-T1-002-158).}% <-this % stops a space
\thanks{Y. Pang and G. Hu are with the School of Electrical and Electronic Engineering, Nanyang
Technological University, 639798, Singapore
        { ypang005@e.ntu.edu.sg, gqhu@ntu.edu.sg}.}%
}

\maketitle
% \thispagestyle{empty}
% \pagestyle{empty}

% \begin{keyword}                           % Five to ten keywords,  
% Distributed optimization; multi-agent system; gradient-free optimization.               % chosen from the IFAC 
% \end{keyword}                             % keyword list or with the 
                                          % help of the Automatica 
                                          % keyword wizard
\begin{abstract}                          % Abstract of not more than 200 words.
This paper presents a special type of distributed optimization problems, where the summation of agents' local cost functions (\textit{i.e.}, global cost function) is convex, but each individual can be non-convex. 
% This problem setting has shown potential benefits in privacy preservation, an important concern commonly raised in networked systems, especially for the differential privacy based objective perturbation methods, where the resulting perturbed function can be non-convex.
Unlike most distributed optimization algorithms by taking the advantages of gradient, the considered problem is allowed to be non-smooth, and the gradient information is unknown to the agents. To solve the problem, a Gaussian-smoothing technique is introduced and a gradient-free method is proposed. 
We prove that each agent's iterate approximately converges to the optimal solution both with probability 1 and in mean, and provide an upper bound on the optimality gap, characterized by the difference between the functional value of the iterate and the optimal value.
The performance of the proposed algorithm is demonstrated by a numerical example and an application in privacy enhancement.
\end{abstract}

% \begin{abstract}
% In this note, we study an online multi-agent optimization problem where the objective functions of agents vary with time. A gradient-free distributed algorithm is proposed
% \end{abstract}


\section{Introduction}
With the prevalence of multi-agent systems, there has been a growing interest in solving the optimization problem in distributed settings recently. One advantage of doing so is that agents access local information and communicate with the neighbors only, making it suitable for the applications with large data size, huge computations and complex network structure. Theoretical works on distributed convex optimization algorithms have also been extensively studied (\textit{e.g.,} see the works\cite{Nedic2015, Shi2015, Xi2016, Liu2017, Hosseini2016, Sun2017}, just name a few). These methods have been shown the effectiveness in many applications, such as parameter estimation and detection , source localization, resource allocation, path-planning \cite{Nowak2003,Ram2010,Lesser2003,Rabbat,Shen2012,DeGennaro2006}, etc.
This paper presents a more general type of distributed convex optimization problems, where the summation of agents' local cost functions (\textit{i.e.}, global cost function) is convex, but each individual can be non-convex. 
This problem setting is motivated from an important concern commonly raised in networked systems -- privacy, especially for the differential privacy based objective perturbation methods\cite{Nozari2018}, where agents exchange some perturbed functions with others before executing distributed algorithms to protect the agents' privacy. The following subsection details this motivation by a simple example.

\subsection{Motivating Example}\label{sec:motivating_ex}
% The privacy has become increasingly important in era of information, especially in the fields of distributed optimization and machine learning. We illustrate how the privacy enhancements motivate this research work by an example inspired from \cite{Abbe2012,Gade2016}.

A typical distributed convex optimization problem among $n$ agents over a graph $\mathcal{G}$ is usually defined as
\begin{align*}
\min F(x) = \sum_{i=1}^n F_i(x),\quad x \in \mathbb{R}^m,
\end{align*}
where the local cost function $F_i$, $i\in\{1,\ldots,n\}$ is assumed to be convex, giving rise to a standard convex setting. 
The graph $\mathcal{G}$ is defined by $\mathcal{G} = \{\mathcal{V},\mathcal{E}\}$, where $\mathcal{V} = \{1, 2, \ldots, n\}$ denotes the set of agents, and $\mathcal{E} \subset \mathcal{V}\times\mathcal{V}$ represents the set of ordered pairs. For any $(i,j)\in\mathcal{E}$, $i,j\in\mathcal{V}$, the information can be transfered from agent $i$ to agent $j$. Assume $(i,i)\in\mathcal{E}$, $\forall i\in\mathcal{V}$. The in-neighbors (respectively, out-neighbors) of agent $i$ are denoted by $\mathcal{N}^{\text{in}}_i = \{j \in \mathcal{V} | (j,i)\in \mathcal{E}\}$ (respectively, $\mathcal{N}^{\text{out}}_i = \{j \in \mathcal{V} | (i,j)\in \mathcal{E}\}$).
To enhance the privacy, each agent randomly generates some (possibly non-convex) functions $F_{i,j}$, $j\in\mathcal{N}_i^{\text{out}}$ and passes them to its out-neighbours. Then, each agent subtracts the self-generated functions from its own local cost function, and combines the functions received from its in-neighbors to form a new local cost function, given by
\begin{align*}
\tilde{F}_i(x) = F_i(x) - \sum_{j\in\mathcal{N}_i^{\text{out}}}F_{i,j}(x) + \sum_{s\in\mathcal{N}_i^{\text{in}}}F_{s,i}(x).
\end{align*}
Obviously, the new local cost function $\tilde{F}_i$ of each agent may be non-convex due to the summation and subtraction of the randomly generated functions by itself and neighbors. However, the new global cost function given by
\begin{align*}
\tilde{F}(x) &= \sum_{i=1}^n \tilde{F}_i(x) = \sum_{i=1}^n\bigg(F_i(x) - \sum_{j\in\mathcal{N}_i^{\text{out}}}F_{i,j}(x) + \sum_{s\in\mathcal{N}_i^{\text{in}}}F_{s,i}(x)\bigg)=\sum_{i=1}^n F_i(x) = F(x)
\end{align*}
is the same as the original global cost function due to the cancellation in the aggregation of local cost functions. Therefore, the new distributed optimization problem given by
\begin{align*}
\min \tilde{F}(x) = \sum_{i=1}^n \tilde{F}_i(x),\quad x \in \mathbb{R}^m
\end{align*}
is equivalent to the original problem. Solving the new distributed optimization problem instead of the original one gives the same optimal solution with the privacy being guaranteed, since all the information provided by each agent in the optimization is purely based on $\tilde{F}_i$ instead of the original $F_i$. Hence, no agents can learn the local cost functions of any other agents based on the provided information. 
% Therefore, to implement the privacy enhance procedures, it is necessary to investigate a special type of distributed optimization problems where the global cost function is convex but local cost functions can be non-convex.

\subsection{Literature Review}

Distributed optimization with non-convex settings have been studied in the literature\cite{Zhu2013,Bianchi2013,Lorenzo2016,Sun2016}. These works considered general non-convex optimization problems where the global cost function, local cost functions and constraints are non-convex. The reported algorithms either obtain in-exact convergence, or need multiple complex steps of implementations. Instead of considering general non-convex settings, it would also be interesting to study a special type of non-convex problems with the global cost function to be convex but each local cost function to be non-convex \cite{Gade2016}, which could potentially enhance privacy. Recent researches on the privacy issue in distributed optimization have been reported in these works \cite{Nozari2018,Lou2018,Huang2015,Hale2015,Han2017,Li2018,Lu2020a}, where almost all of them employ randomized perturbation techniques in the message sharing steps \cite{Huang2015,Hale2015,Han2017,Li2018} or in the cost functions \cite{Nozari2018}. 
It should be noted that the perturbation process in the cost function may affect the smoothness and convexity of the local cost fuction, which makes it non-differentiable and non-convex. Therefore, it motivates us to study a special type of non-convex problems with a convex global cost but possibly non-convex local costs, and propose a distributed algorithm to solve it. 
We remark that the applications of the typical distributed convex optimization are still feasible to our proposed algorithm. Moreover, as an extension to the typical distributed convex optimization, our proposed algorithm can be applied to the cases where the local cost functions are non-convex. Though being limited, one direct application could be about the privacy protection in networked systems, as discussed in Sec.~\ref{sec:motivating_ex}. Other applications could be related to the distributed non-convex problems where convexification is appropriate. In these applications, only the global cost function is needed to be convexified in our work while all local cost functions are needed to be convexified in the typical distributed convex optimization.

From the perspective of algorithmic development, it is noted that all the aforementioned algorithms considered the problem where the gradient information of the cost functions is directly accessible, and hence proposed a variety of gradient-based methods, which have been shown to be effective, however, cannot be applied if the gradient information is not available. This motivates us to study the gradient-free optimization.
The idea of the gradient-free approach was initially brought out in the work \cite{Matyas1965}, and received a renewed attention recently \cite{Nesterov2017,Duchi2015,Yuan2015,Li2015,Yuan2015a,Chen2017,Pang2017,Pang2020,Pang2018,Sahu2019,Sahu2018}. Specifically, the work \cite{Nesterov2017} considered a general convex optimization problem, and proposed a Gaussian smoothing technique, where a randomized gradient-free oracle was constructed based on two-point evaluations. Convergence results have been derived for the cost function with different degrees of smoothness. The work \cite{Duchi2015} characterized the optimal rates of convergence based on multiple noisy function evaluations, and presented simple randomization-based algorithms to achieve these optimal rates. For distributed optimization problems, where the gradient-free oracle was integrated with the subgradient method \cite{Yuan2015,Li2015} and the push-sum algorithm \cite{Yuan2015a}, and was further extended to a two-sided gradient-free oracle \cite{Chen2017}. Our recent work \cite{Pang2017} applied the gradient-free technique in an unconstrained distributed optimization problem but using less stringent weight matrix, which was further studied in constrained scenarios \cite{Pang2020,Pang2018}. The works \cite{Sahu2019,Sahu2018} also studied gradient-free distributed optimization for strongly convex cost functions based on an adaptive probabilistic sparsifying communications protocol. In addition to the distributed optimization, vast results on these gradient-free techniques have been reported in bandit online optimization literature\cite{Pang2019b,Li2020,Zhang2020,Yi2020,Lu2020}. 
It should be noted that all the aforementioned works on gradient-free algorithms mainly focus on the problems where the cost function of each agent is (strongly) convex. Little attention has been received for the problems with non-convex local cost functions.

In this paper, we propose a gradient-free algorithm using a Gaussian smoothing technique \cite{Nesterov2017} to solve a special type of non-convex optimization problems with a convex global cost but possibly non-convex local costs.
The major contributions as compared to the literature are threefold.
\begin{enumerate}
\item From the perspective of the problem setup, unlike most existing distributed optimization algorithms with the assumption of convexity on the cost functions, the proposed approach only requires the convexity of the global cost function, but each local cost can be non-convex. It should be noted that the relaxation of the convexity of cost functions are non-trivial. Practically, non-convex settings are more closely related to the applications in the real world. Theoretically, it is more challenging to prove the optimality without the assumption of convexity.
\item From the perspective of the algorithm development, the proposed algorithm does not require the knowledge of first-order information (\textit{i.e.}, gradient or subgradient information). Only zero-order information (\textit{i.e.}, value of the local cost function) is needed throughout the optimization process. Hence, it can be applied to those problems where the gradient information is difficult to compute or even does not exist, which increases the range of the applications. Besides, the proposed algorithm adopted the surplus-based method\cite{Cai2012}, which removes the typical doubly stochastic requirement on the weight matrix. By doing so, it further extends the feasible range of the algorithm to those directed graphs which admit no corresponding doubly stochastic weight matrices\cite{Gharesifard2010}.
\item From the perspective of the achieved results, in addition to the convergence of the functional value of each agent's estimate in the mean sense\cite{Yuan2015,Li2015,Yuan2015a,Pang2017,Pang2020,Pang2018}, we also obtain the convergence of each agent's estimate both with probability 1 and in mean, which is more comprehensive. Theoretical analysis is provided to show that the poposed algorithm can achieve the convergence to an approximate optimal solution, where the optimality
gap is characterized by the smoothing parameter. The performance of the algorithm is demonstrated by a numerical example and an application in privacy enhancement, where influences of different smoothing parameters, number of agents, and problem dimensions on the convergence results are investigated.
\end{enumerate}

In the following sections, the problem is defined in Sec.~\ref{sec:problem_formulation}. Main results of the algorithm and established properties are reported in Sec.~\ref{sec:distr_opt}. In Sec.~\ref{sec:simulation}, the performance of the proposed approach is justified by a numerical example followed by an application in privacy enhancement.

\textbf{Notations:}
We denote the set of real numbers and $m$-dimensional column vectors by $\mathbb{R}$ and $\mathbb{R}^m$, respectively, the element in the $i$-th row and $j$-th column of a matrix $W$ by $[W]_{ij}$, the derivative of a differentiable function $f(x)$ by $\nabla f(x)$, the standard Euclidean norm of a vector $x$ by $\|x\|$, and the number of elements in a set $\mathcal{V}$ by $|\mathcal{V}|$.

\section{Problem Formulation}\label{sec:problem_formulation}

In this section, the problem is formally defined, followed by some preliminary results.

\subsection{Problem Definition}

We consider the following distributed optimization problem among $n$ agents over the directed graph $\mathcal{G}$:  
\begin{align}\label{eq:cost_function}
  \min f(x) = \sum_{i=1}^n f_i(x),
\end{align}
where $x \in \mathbb{R}^m$ is a global decision vector, and $f_i$ is a local cost function. 
We denote the (non-empty) optimal solution set of problem \eqref{eq:cost_function} by $\mathcal{X}^\star= \arg\min_{x\in \mathbb{R}^m}f(x)$, and the optimal value by $f^\star$, \textit{i.e.}, $f^\star = f(x^\star)$ for any $x^\star\in\mathcal{X}^\star$.

The following standard assumptions are made throughout the paper:
\begin{Assumption}\label{assumption_graph}
The directed graph $\mathcal{G}$ is strongly connected.
% , \textit{i.e.,} there is a path in each direction between each pair of vertices of the graph.
\end{Assumption}
\begin{Assumption}\label{assumption_local_f_lipschitz}
The global cost function $f$ is convex. Each local cost function $f_i$ can be non-convex, and is assumed to be Lipschitz continuous with a constant $\hat{D}$,
\textit{i.e.,} $\forall x, y \in \mathbb{R}^m$, there exists a constant $\hat{D}>0$ such that $|f_i(x) - f_i(y)|\leq \hat{D}\|x-y\|$, $\forall i\in\mathcal{V}$.
\end{Assumption}

\subsection{Preliminaries}
Motivated by the Gaussian smoothing technique \cite{Nesterov2017}, a smoothed version problem of (\ref{eq:cost_function}) can be defined as
\begin{align}\label{eq:smooth_cost_function}
  \min f_\mu(x) = \sum_{i=1}^n f_{i,\mu}(x),\quad x \in \mathbb{R}^m,
\end{align}
where $f_{i,\mu}(x)$ is a Gaussian smoothed function of $f_i(x)$, given by
\begin{align*}
  f_{i,\mu}(x) = \frac1{\kappa}\int_{\mathbb{R}^m} f_i(x+\mu\zeta)\exp\bigg({-\frac12\|\zeta\|^2}\bigg)d\zeta,
\end{align*}
with $\kappa = \int_{\mathbb{R}^m} \exp({-\frac12\|\zeta\|^2})d\zeta = (2\pi)^{m/2}$ and $\mu \geq 0$ is a smoothing parameter. Likewise, we denote the (non-empty) optimal solution set of problem \eqref{eq:smooth_cost_function} by $\mathcal{X}^\star_\mu= \arg\min_{x\in \mathbb{R}^m}f_\mu(x)$, and the optimal value by $f^\star_\mu$, \textit{i.e.}, $f^\star_\mu = f_\mu(x^\star_\mu)$ for any $x^\star_\mu\in\mathcal{X}^\star_\mu$. Then, the randomized gradient-free oracle of $f_i(x)$ can be designed as follows\cite{Nesterov2017}:
\begin{align*}
g^i_{\mu}(x) = \frac{f_i(x+\mu\zeta^i)-f_i(x)}{\mu}\zeta^i,
\end{align*}
where $\zeta^i \in \mathbb{R}^m$ is a normally distributed Gaussian vector.
The properties of the functions $g^i_{\mu}(x)$ and $f_{i,\mu}(x)$ are presented in the following two lemmas:
\begin{Lemma}\label{lemma:property_f_mu}
(see Sec.~2\cite{Nesterov2017}) Under Assumption~\ref{assumption_local_f_lipschitz}, we have that
\begin{enumerate}
\item function $f_{i,\mu}$ is not necessarily convex due to non-convex $f_i$, but $f_\mu = \sum_{i=1}^n f_{i,\mu}$ is convex due to convex $f$. Moreover, $f_{i,\mu}(x)$ satisfies $|f_{i,\mu}(x) - f_i(x)| \leq \sqrt{m}{\mu}\hat{D}$,
\item for $\mu>0$, function $f_{i,\mu}(x)$ is differentiable, and Lipschitz continuous with a constant $\hat{D}_\mu\leq\hat{D}$. Its gradient $\nabla f_{i,\mu}(x)$ is Lipschitz continuous with a constant $\hat{L} = \frac{\sqrt{m}\hat{D}}{\mu}$, \textit{i.e.}, $\|\nabla f_{i,\mu}(x) - \nabla f_{i,\mu}(y)\| \leq \hat{L}\|x - y\|$, and holds that $\nabla f_{i,\mu}(x) = \mathbf{E}[g^i_{\mu}(x)]$,
\item the randomized gradient-free oracle $g^i_{\mu}(x)$ holds that $\mathbf{E}[\|g^i_{\mu}(x)\|] \leq \sqrt{\mathbf{E}[\|g^i_{\mu}(x)\|^2]} \leq \mathcal{B}$, where $\mathcal{B} = (m+4)\hat{D}$.
\end{enumerate}
\end{Lemma}

\section{Main Results}\label{sec:distr_opt}

% In this section, we will review our previously developed algorithm -- randomized gradient-free directed-distributed gradient descent (RGF-D-DGD) in \cite{Pang2017} for completeness, followed by the convergence analysis.
In this section, the proposed method is introduced in details. Then, we rigorously analyze the convergence property.

\subsection{Randomized Gradient-Free DGD Method} 
The details of our proposed algorithm are described as follows.

Every agent $j$ passes the information of its decision $x^j_k$ with a weighted information $[W_c]_{ij}y^j_k$ to its out-neighbor $i \in \mathcal{N}^\text{out}_j$ at time $k$, where $y^j$ is an auxiliary variable of agent $j$ to offset the shift caused by the unbalanced (non-doubly stochastic) weighting structure. Upon receiving the information, every agent $i$ proceeds to update its decision $x_{k+1}^i$ and variable $y_{k+1}^i$ based on the following updating laws:
\begin{subequations}\label{eq:update_law}
\begin{align}
x_{k+1}^i &= \sum_{j=1}^n [W_r]_{ij}x^j_k + \epsilon y^i_k - \alpha_k g^i_{\mu}(x^i_k),\label{eq:update_x}\\
y_{k+1}^i &= x^i_k - \sum_{j=1}^n [W_r]_{ij}x^j_k +  \sum_{j=1}^n [W_c]_{ij} y^j_k - \epsilon y^i_k,\label{eq:update_y}
\end{align}
\end{subequations}
where $g^i_{\mu}(x^i_k)$ is the randomized gradient-free oracle
\begin{align}
g^i_{\mu}(x^i_k) = \frac{f_i(x^i_k+\mu\zeta^i_k)-f_i(x^i_k)}{\mu}\zeta^i_k, \label{grad_oracle}
\end{align}
matrix $W_r$ is row-stochastic ($\sum_{j=1}^n[W_r]_{ij} = 1$ for all $i \in \mathcal{V}$), and matrix $W_c$ is column-stochastic ($\sum_{i=1}^n[W_c]_{ij} = 1$ for all $j \in \mathcal{V}$); $\epsilon$ is a small positive number; and $\alpha_k \geq 0$ is a step-size satisfying
\begin{equation}\label{eq:step-size_condition}
\sum_{k=0}^\infty \alpha_k = \infty, \quad \sum_{k=0}^\infty \alpha^2_k < \infty.
\end{equation}

\begin{Remark}
In the case where the global decision variable is constrained by a general convex set, the proposed algorithm can be extended to solve the problem by introducing a projection operator in \eqref{eq:update_x}\cite{Xi2016,Pang2020}.
\end{Remark}

\subsection{Convergence Analysis}

The convergence properties of our proposed algorithm are detailed in this part.
Denote by $\mathcal{F}_k$ the $\sigma$-algebra generated by $\{\zeta^i_0,\zeta^i_1,\ldots,\zeta^i_{k-1},i\in\mathcal{V}\}$, $k\geq1$.
We write (\ref{eq:update_x}) and (\ref{eq:update_y}) in a compact form as
\begin{align}\label{eq:d-dgd}
\theta^i_{k+1} = \sum_{j=1}^{2n}[W]_{ij}\theta^j_k - \alpha_k g^i_k,
\end{align}
where $\theta^i_k = x^i_k$ for $1 \leq i \leq n$, $\theta^i_k = y^{i-n}_k$ for $n+1 \leq i \leq 2n$; $g^i_k = g^i_{\mu}(x^i_k)$ for $1 \leq i \leq n$, $g^i_k = \mathbf{0}_n$ for $n+1 \leq i \leq 2n$; and $W = [\begin{smallmatrix} W_r & \epsilon I \\ I-W_r & W_c - \epsilon I\end{smallmatrix}]$. The following lemma presents a convergence property on matrix $W$.

\begin{Lemma}\label{lemma:A_matrix}
(Lemma 1\cite{Xi2016}) Suppose Assumption~\ref{assumption_graph} holds. Let $W$ be the weight matrix in \eqref{eq:d-dgd} with $\epsilon$ satisfying $\epsilon \in (0,\bar{\epsilon})$, where $\bar{\epsilon} = \frac1{(20+8n)^n}(1-|\lambda_3|)^n$, $\lambda_3$ is the third largest eigenvalue of $W$ by setting $\epsilon =0$. Then $\forall i,j \in \{1,\ldots,2n\}$, the entry $[W^k]_{ij}$ converges at a geometric rate. Specifically, for $k\geq 1$,
\begin{equation*}
\left\|W^k - \begin{bmatrix} \frac{\mathbf{1}_n\mathbf{1}^T_n}n &\frac{\mathbf{1}_n\mathbf{1}^T_n}n \\ \mathbf{0} & \mathbf{0} \end{bmatrix}\right\|_\infty \leq \Gamma \gamma^k,
\end{equation*}
where $\Gamma > 0$ and $0 < \gamma < 1$ are some constants.
\end{Lemma}

Now, the average of all agents' information at time $k$ can be defined by
\begin{equation}\label{eq:def_z_bar}
\bar{\theta}_k = \frac1n\sum_{i=1}^{2n}\theta^i_k = \frac1n\sum_{i=1}^nx^i_k + \frac1n\sum_{i=1}^ny^i_k.
\end{equation}

In the following theorem, we would like to characterize the consensus property of the algorithm: the state information of all agents $x^i_k, i \in \mathcal{V}$ will converge to their average $\bar{\theta}_k$ with probability 1. 

\begin{Theorem}\label{theorem:consensus}
Under Assumptions~\ref{assumption_graph} and \ref{assumption_local_f_lipschitz},  sequence $\{x^i_k\}_{k\geq0}$, $i \in \mathcal{V}$ is obtained from the update law \eqref{eq:update_law} with $\epsilon$ in $W$ satisfying $\epsilon \in (0,\bar{\epsilon})$, where $\bar{\epsilon} = \frac1{(20+8n)^n}(1-|\lambda_3|)^n$, $\lambda_3$ is the third largest eigenvalue of $W$ by setting $\epsilon =0$, $\Gamma > 0$ and $0 < \gamma < 1$ are some constants, and the step-size sequence $\{\alpha_k\}_{k\geq0}$ satisfying \eqref{eq:step-size_condition}. Then, we have\footnote{We use `a.s.' for `almost surely'. A sequence of random vectors $\{v_k\}_{k\geq0}$ converges to $v$, if the probability of $\lim_{k\to\infty}v_k = v$ is 1.}
\begin{enumerate}
\item $\begin{aligned}[t] \sum_{k=1}^\infty\alpha_k\sum_{i=1}^n\|x^i_k - \bar{\theta}_k\|<\infty\end{aligned}$ a.s.;
\item the sequence $\{x^i_k-\bar{\theta}_k\}_{k\geq0}$ converges to 0 with probability 1 and in mean, 
\end{enumerate}
where $\bar{\theta}_k$ is defined in \eqref{eq:def_z_bar}.
\end{Theorem}
\begin{Proof}
Recursively expanding \eqref{eq:d-dgd} for $1\leq i \leq 2n$, $k\geq1$ and noting that $g^i_k = \mathbf{0}_n$ for $n+1\leq i\leq 2n$, we have
\begin{align*}
\theta^i_k = \sum_{j=1}^{2n}[W^k]_{ij}\theta^j_0 - \sum_{r=1}^{k-1}\sum_{j=1}^n[W^{k-r}]_{ij}\alpha_{r-1}g^j_{r-1}-\alpha_{k-1} g^i_{k-1},
\end{align*}
which by definition \eqref{eq:def_z_bar} gives
\begin{align*}
\bar{\theta}_k = \frac1n\sum_{j=1}^{2n}\theta^j_0 - \frac1n\sum_{r=1}^{k-1}\sum_{j=1}^n\alpha_{r-1}g^j_{r-1}-\frac{\alpha_{k-1}}n\sum_{j=1}^n g^j_{k-1},
\end{align*}
where we have applied $\sum_{i=1}^{2n}[W^k]_{ij} = 1$ for any $1\leq j \leq 2n$, $k\geq1$. Taking the subtraction for the above two relations, we have for $i\in\mathcal{V}$, $k\geq1$
\begin{align}
\|x^i_k-\bar{\theta}_k\| &\leq \sum_{j=1}^{2n}\bigg\|[W^k]_{ij}-\frac1n\bigg\|\theta^j_0 + \sum_{r=1}^{k-1}\sum_{j=1}^n\bigg\|[W^{k-r}]_{ij}-\frac1n\bigg\|\alpha_{r-1}\|g^j_{r-1}\|\nonumber\\
&\quad\quad+\frac{(n-1)\alpha_{k-1}}n \|g^i_{k-1}\|+\frac{\alpha_{k-1}}n\sum_{j\neq i} \|g^j_{k-1}\|\nonumber\\
&\leq 2n\delta\Gamma\gamma^k + n\Gamma\sum_{r=1}^{k-1}\gamma^{k-r}\alpha_{r-1}\|g^j_{r-1}\|+\alpha_{k-1}\sum_{j=1}^n\|g^j_{k-1}\|, \label{eq:E_z-E_b}
\end{align}
where $\delta=\max\{\|x^j_0\|,\|y^j_0\|,j\in\mathcal{V}\}$ and the second inequality is due to Lemma~\ref{lemma:A_matrix}. Applying Lemma~\ref{lemma:property_f_mu}-(3), it follows from \eqref{eq:E_z-E_b} that
\begin{align*}
\alpha_k\sum_{i=1}^n\mathbf{E}[\|x^i_k - \bar{\theta}_k\|] 
\leq 2n^2\delta\Gamma\alpha_k\gamma^k + n^2\mathcal{B}\Gamma\sum_{r=1}^{k-1}\gamma^{k-r}\alpha_k\alpha_{r-1}+ n^2\mathcal{B}\alpha_k\alpha_{k-1}.
\end{align*}
Summing over $k = 1, \ldots, t$ and recalling that 
\begin{align*}
\sum_{k=1}^t \alpha_k \gamma^k &\leq \frac12\sum_{k=1}^t\alpha^2_k+\frac{\gamma^2}{2(1-\gamma^2)},
\sum_{k=1}^t\sum_{r=1}^{k-1}\gamma^{k-r}\alpha_k\alpha_{r-1} \leq \frac{\gamma}{1-\gamma}\sum_{k=1}^t\alpha^2_k,
\sum_{k=1}^t\alpha_k\alpha_{k-1} \leq \sum_{k=0}^t\alpha^2_k,
\end{align*}
from Lemma~3\cite{Pang2017}, we have
\begin{align}
\sum_{k=1}^t\alpha_k\sum_{i=1}^n\mathbf{E}[\|x^i_k - \bar{\theta}_k\|] \leq n^2\delta\Gamma\bigg(\sum_{k=1}^t\alpha^2_k+\frac{\gamma^2}{1-\gamma^2}\bigg)+ \frac{n^2\mathcal{B}\Gamma\gamma}{1-\gamma}\sum_{k=1}^t\alpha^2_k+ n^2\mathcal{B}\sum_{k=0}^t\alpha^2_k.\label{eq:eta_k_sum}
\end{align}
Letting $t\to \infty$, and noting that $\sum_{k=0}^\infty\alpha_k^2<\infty$, it follows from \eqref{eq:eta_k_sum} that $\sum_{k=1}^\infty\alpha_k\sum_{i=1}^n\mathbf{E}[\|x^i_k - \bar{\theta}_k\|]<\infty$, which by the monotone convergence theorem implies $\mathbf{E}[\sum_{k=1}^\infty\alpha_k\sum_{i=1}^n\|x^i_k - \bar{\theta}_k\|]<\infty$. Therefore, we conclude $\sum_{k=1}^\infty\alpha_k\sum_{i=1}^n\|x^i_k - \bar{\theta}_k\|<\infty$ a.s, which completes the proof of part (1).

For part (2), taking the square on both sides of \eqref{eq:E_z-E_b}, and applying $\|a+b+c\|^2\leq3(\|a\|^2+\|b\|^2+\|c\|^2)$, we have for $i\in\mathcal{V}$, $k\geq1$
\begin{align*}
\|x^i_k-\bar{\theta}_k\|^2 &\leq 12n^2\delta^2\Gamma^2\gamma^{2k} + 3n^2\Gamma^2\bigg(\sum_{r=1}^{k-1}\gamma^{k-r}\alpha_{r-1}\|g^j_{r-1}\|\bigg)^2+3n\alpha^2_{k-1}\sum_{j=1}^n\|g^j_{k-1}\|^2,
\end{align*}
where the second term follows that
\begin{align*}
\bigg(\sum_{r=1}^{k-1}\gamma^{k-r}\alpha_{r-1}\|g^j_{r-1}\|\bigg)^2 \leq \frac\gamma{1-\gamma}\sum_{r=1}^{k-1}\gamma^{k-r}\alpha^2_{r-1}\|g^j_{r-1}\|^2
\end{align*}
based on Cauchy-Schwarz inequality. Taking the total expectation and summing over from $k=1$ to $\infty$, we obtain
\begin{align*}
\sum_{k=1}^\infty\mathbf{E}[\|x^i_k-\bar{\theta}_k\|^2] &\leq \frac{12n^2\delta^2\Gamma^2\gamma^2}{1-\gamma^2} + \frac{3n^2\mathcal{B}^2\Gamma^2\gamma}{1-\gamma}\sum_{k=1}^\infty\sum_{r=1}^{k-1}\gamma^{k-r}\alpha^2_{r-1}+3n^2\mathcal{B}^2\sum_{k=1}^\infty\alpha^2_{k-1}\\
&\leq \frac{12n^2\delta^2\Gamma^2\gamma^2}{1-\gamma^2} + \bigg(\frac{3n^2\mathcal{B}^2\Gamma^2\gamma^2}{(1-\gamma)^2}+3n^2\mathcal{B}^2\bigg)\sum_{k=0}^\infty\alpha^2_{k-1},
\end{align*}
where we have applied $\sum_{k=1}^\infty\sum_{r=1}^{k-1}\gamma^{k-r}\alpha^2_{r-1}\leq\frac\gamma{1-\gamma}\sum_{k=0}^\infty\alpha^2_k$. Since $\sum_{k=0}^\infty\alpha_k^2<\infty$, we have
\begin{align}
\sum_{k=1}^\infty\mathbf{E}[\|x^i_k-\bar{\theta}_k\|^2]<\infty,\label{eq:E_z-E_b_sq}
\end{align}
which by the monotone convergence theorem implies $\mathbf{E}[\sum_{k=1}^\infty\|x^i_k-\bar{\theta}_k\|^2]<\infty$. Therefore, we conclude the sequence $\{x^i_k-\bar{\theta}_k\}_{k\geq0}$ converges to 0 with probability 1 and in mean for all $i \in \mathcal{V}$.
\end{Proof}

% \begin{Remark}
% Theorem~\ref{theorem:consensus} is a characterization of the consensus property of the algorithm, which implies that the state information of all agents $x^i_k, i \in \mathcal{V}$ will converge to their average $\bar{\theta}_k$ with probability 1. 
% \end{Remark}

The following Lemma is the key to the convergence of the proposed algorithm, where the relation between two successive updates is developed.
\begin{Lemma}\label{lemma:successive_update}
Under Assumptions~\ref{assumption_graph} and \ref{assumption_local_f_lipschitz}, sequence $\{\theta^i_k\}_{k\geq0}$ is obtained from the update law (\ref{eq:d-dgd}). Let the step-size sequence $\{\alpha_k\}_{k\geq0}$ be chosen such that condition (\ref{eq:step-size_condition}) is satisfied. Then, for $k \geq 0$, and any $\theta^\star\in\mathcal{X}^\star_\mu$, it holds that
\begin{align}
\mathbf{E}[\|\bar{\theta}_{k+1}-\theta^\star\|^2|\mathcal{F}_k] &\leq - \frac{2\alpha_k}n(f_\mu(\bar{\theta}_k) - f^\star_\mu)+\bigg(1+\frac{2\hat{L}\alpha_k}n\sum_{i=1}^n\|x^i_k - \bar{\theta}_k\|\bigg)\|\bar{\theta}_k-\theta^\star\|^2\nonumber\\
&\quad+\frac{2\hat{L}\alpha_k}n\sum_{i=1}^n\|x^i_k - \bar{\theta}_k\|+ \mathcal{B}^2\alpha^2_k,\label{eq:successive_update}
\end{align}
where $\mathcal{B}$ and $\hat{L}$ are defined in Lemma~\ref{lemma:property_f_mu}.
\end{Lemma}
\begin{Proof}
By the definition of $\bar{\theta}_k$, it can be obtained that
\begin{align}
\bar{\theta}_{k+1} = \bar{\theta}_k- \frac{\alpha_k}n\sum_{i=1}^ng^i_{\mu}(x^i_k). \label{eq:zbar_iteration}
\end{align}
For any $\theta^\star\in\mathcal{X}^\star_\mu$, subtracting both sides by $\theta^\star$ and taking the norm, it yields that
\begin{align*}
\|\bar{\theta}_{k+1}-\theta^\star\|^2 &= \|\bar{\theta}_k-\theta^\star\|^2 + \bigg\|\frac{\alpha_k}n\sum_{i=1}^ng^i_{\mu}(x^i_k)\bigg\|^2- \frac{2\alpha_k}n\sum_{i=1}^ng^i_{\mu}(x^i_k)^T(\bar{\theta}_k-\theta^\star).
\end{align*}
Thus, we have
\begin{align}
\mathbf{E}[\|\bar{\theta}_{k+1}-\theta^\star\|^2|\mathcal{F}_k] 
\leq& \|\bar{\theta}_k-\theta^\star\|^2 + \frac{\alpha^2_k}n\sum_{i=1}^n\mathbf{E}[\|g^i_{\mu}(x^i_k)\|^2|\mathcal{F}_k] -\frac{2\alpha_k}n\sum_{i=1}^n\mathbf{E}[g^i_{\mu}(x^i_k)|\mathcal{F}_k]^T(\bar{\theta}_k-\theta^\star)\nonumber\\
\leq& \|\bar{\theta}_k-\theta^\star\|^2 + \mathcal{B}^2\alpha^2_k- \frac{2\alpha_k}n\sum_{i=1}^n\nabla f_{i,\mu}(x^i_k)^T(\bar{\theta}_k-\theta^\star),\label{eq:z_bar_minus_z_star}
\end{align}
where we have applied the results from Lemma~\ref{lemma:property_f_mu}. For the last term
\begin{align}
&\sum_{i=1}^n\nabla f_{i,\mu}(x^i_k)^T(\bar{\theta}_k-\theta^\star) = \sum_{i=1}^n\nabla f_{i,\mu}(\bar{\theta}_k)^T(\bar{\theta}_k-\theta^\star) + \sum_{i=1}^n\big(\nabla f_{i,\mu}(x^i_k)-\nabla f_{i,\mu}(\bar{\theta}_k)\big)^T(\bar{\theta}_k-\theta^\star).\label{eq:grad_times_z_bar_minus_z_star}
\end{align}
Recalling that $f_{\mu}$ is convex due to Lemma~\ref{lemma:property_f_mu}-(1), we have the first term of \eqref{eq:grad_times_z_bar_minus_z_star}
\begin{align*}
\sum_{i=1}^n\nabla f_{i,\mu}(\bar{\theta}_k)^T(\bar{\theta}_k-\theta^\star) &=\nabla f_{\mu}(\bar{\theta}_k)^T(\bar{\theta}_k-\theta^\star) \geq f_\mu(\bar{\theta}_k) - f_\mu^\star,
\end{align*}
and the second term of \eqref{eq:grad_times_z_bar_minus_z_star}
\begin{align*}
\sum_{i=1}^n\big(\nabla f_{i,\mu}(x^i_k)-\nabla f_{i,\mu}(\bar{\theta}_k)\big)^T(\bar{\theta}_k-\theta^\star)
\geq&- \hat{L}\sum_{i=1}^n\|x^i_k - \bar{\theta}_k\|\|\bar{\theta}_k-\theta^\star\|\\
\geq&- \hat{L}\sum_{i=1}^n\|x^i_k - \bar{\theta}_k\|(1+\|\bar{\theta}_k-\theta^\star\|^2),
\end{align*}
where we have used Lemma~\ref{lemma:property_f_mu}-(2) and Young's inequality ($\|\mathbf{z}\|\leq2\|\mathbf{z}\|\leq1+\|\mathbf{z}\|^2$). Combining the above results, it follows from (\ref{eq:grad_times_z_bar_minus_z_star}) that
\begin{align}
&\sum_{i=1}^n\nabla f_{i,\mu}(x^i_k)^T(\bar{\theta}_k-\theta^\star) \geq f_\mu(\bar{\theta}_k) - f^\star_\mu - \hat{L}\sum_{i=1}^n\|x^i_k - \bar{\theta}_k\|-\hat{L}\sum_{i=1}^n\|x^i_k - \bar{\theta}_k\|\|\bar{\theta}_k-\theta^\star\|^2.\label{eq:grad_times_z_bar_minus_z_star_2}
\end{align}
The conclusion follows by substituting (\ref{eq:grad_times_z_bar_minus_z_star_2}) into (\ref{eq:z_bar_minus_z_star}).
\end{Proof}

Now, we are ready to establish the optimality property, as stated in the theorem below.
\begin{Theorem}\label{theorem:optimality}
Under Assumptions~\ref{assumption_graph} and \ref{assumption_local_f_lipschitz}, sequence $\{x^i_k\}_{k\geq0}$ is obtained from the update law (\ref{eq:d-dgd}) with $\epsilon$ in $W$ satisfying $\epsilon \leq \min(\bar{\epsilon},\frac{1-\gamma}{2n\Gamma\gamma})$. Let the step-size sequence $\{\alpha_k\}_{k\geq0}$ be chosen such that condition (\ref{eq:step-size_condition}) is satisfied. Then, for $i\in\mathcal{V}$, we have
\begin{enumerate}
\item $\begin{aligned}[t]\{x^i_k\}_{k\geq0}\end{aligned}$ converges to some $x^\star_\mu\in\mathcal{X}^\star_\mu$ with probability 1 and in mean;
\item $\begin{aligned}[t]\lim_{k\to\infty}\mathbf{E}[f(x^i_k)] - f^\star \leq 2n\sqrt{m}{\mu}\hat{D}\end{aligned}$.
\end{enumerate}
\end{Theorem}
\begin{Proof}
We rewrite (\ref{eq:successive_update}) in Lemma~\ref{lemma:successive_update} as
\begin{equation}\label{eq:f_z_bar_minus_f_mu}
\mathbf{E}[{v}_{k+1}|\mathcal{F}_k] \leq (1+\eta_k){v}_k - {b}_k + {c}_k,
\end{equation}
where 
\begin{align*}
{v}_k &= \|\bar{\theta}_k-\theta^\star\|^2,\quad
{b}_k = \frac{2\alpha_k}n(f_\mu(\bar{\theta}_k) - f^\star_\mu),\\
{c}_k &= \frac{2\hat{L}\alpha_k}n\sum_{i=1}^n\|x^i_k - \bar{\theta}_k\|+ \mathcal{B}^2\alpha^2_k,\quad
\eta_k = \frac{2\hat{L}\alpha_k}n\sum_{i=1}^n\|x^i_k - \bar{\theta}_k\|.
\end{align*}
Based on Theorem~\ref{theorem:consensus}-(1) and $\sum_{k=1}^\infty\alpha^2_k<\infty$, we can show that
\begin{align*}
&\sum_{k=1}^\infty {c}_k= \frac{2\hat{L}}n\sum_{k=1}^\infty\alpha_k\sum_{i=1}^n\|x^i_k - \bar{\theta}_k\|+ \mathcal{B}^2\sum_{k=1}^\infty\alpha^2_k<\infty, \text{ a.s.},\\
&\sum_{k=1}^\infty \eta_k = \frac{2\hat{L}}n\sum_{k=1}^\infty\alpha_k\sum_{i=1}^n\|x^i_k - \bar{\theta}_k\| <\infty, \text{ a.s.}
\end{align*}
Then, we introduce the following lemma to facilitate the proof of convergence.
% (\textit{e.g.,} \cite{Gade2016, Xi2016a})

\begin{Lemma}\label{lemma:sequence_relation}
(Robbins-Siegmund's Lemma \cite{Robbins1985,Polyak1987}) Let ${v}_k, {b}_k, {c}_k,\eta_k$ be non-negative random variables satisfying that 
$\sum_{k=0}^\infty \eta_k < \infty$ a.s., $\sum_{k=0}^\infty {c}_k < \infty$ a.s., and $\mathbf{E}[{v}_{k+1}|\mathcal{F}_k] \leq (1+\eta_k){v}_k - {b}_k + {c}_k$ a.s., 
where $\mathbf{E}[{v}_{k+1}|\mathcal{F}_k]$ denotes the conditional expectation for the given ${v}_0,\ldots,{v}_k$, ${b}_0,\ldots,{b}_k$, ${c}_0,\ldots,{c}_k$, $\eta_0,\ldots,\eta_k$. Then,
\begin{enumerate}
\item $\begin{aligned}[t]\{{v}_k\}_{k\geq0}\end{aligned}$ converges a.s., and $\sup_{k\geq0}\mathbf{E}[{v}_k]<\infty$;
\item $\begin{aligned}[t]\sum_{k=0}^\infty {b}_k < \infty\end{aligned}$ a.s.
\end{enumerate}
\end{Lemma}

\begin{Remark}
The conclusions that $\{{v}_k\}_{k\geq0}$ converges a.s. and $\sum_{k=0}^\infty {b}_k < \infty$ a.s. are well-known in Robbins-Siegmund's Lemma \cite{Robbins1985,Polyak1987}, while the proof of $\sup_{k\geq0}\mathbf{E}[{v}_k]<\infty$ can be referred to Theorem~2.3.5\cite{Gadat2018}.
\end{Remark}

In light of Lemma~\ref{lemma:sequence_relation}, we obtain that
\begin{subequations}
\begin{align}
\{\|\bar{\theta}_k-\theta^\star\|^2\}_{k\geq0} \text{ converges a.s., and } \sup_{k\geq0}\mathbf{E}[\|\bar{\theta}_k-\theta^\star\|^2]<\infty, \text{ for any } \theta^\star\in\mathcal{X}^\star_\mu, \label{eq:result1}\\
\sum_{k=1}^\infty \alpha_k(f_\mu(\bar{\theta}_k) - f^\star_\mu)< \infty \text{ a.s}. \label{eq:result2}
\end{align}
\end{subequations}
Since the step-size sequence satisfies $\sum_{k=0}^\infty \alpha_k = \infty$ and $f_\mu(\bar{\theta}_k) \geq f^\star_\mu$, it follows from \eqref{eq:result2} that $\liminf_{k\to\infty}f_\mu(\bar{\theta}_k) = f_\mu^\star$ a.s. Let $\{\bar{\theta}_\ell\}_{\ell\geq0}$ be a subsequence of $\{\bar{\theta}_k\}_{k\geq0}$ such that
\begin{align}
\lim_{\ell\to\infty}f_\mu(\bar{\theta}_\ell)=\liminf_{k\to\infty}f_\mu(\bar{\theta}_k) = f_\mu^\star \text{ a.s.} \label{eq:result3}
\end{align}
It follows from \eqref{eq:result1} that $\{\bar{\theta}_k\}_{k\geq0}$ is bounded a.s. Without loss of generality, we may assume $\{\bar{\theta}_\ell\}_{\ell\geq0}$ converges a.s. to some $\tilde{x}$ (otherwise, we may choose a subsequence of $\{\bar{\theta}_\ell\}_{\ell\geq0}$ such that it converges a.s.). By continuity of $f_\mu$, we have $f(\bar{\theta}_\ell)$ converges a.s. to $f(\tilde{x})$. In view of \eqref{eq:result3}, we further obtain $f(\tilde{x})=f^\star_\mu$, \textit{i.e.}, $\tilde{x}\in\mathcal{X}^\star_\mu$. Letting $\theta^\star = \tilde{x}$ in \eqref{eq:result1} and considering the sequence $\{\|\bar{\theta}_k-\tilde{x}\|^2\}_{k\geq0}$, it has a subsequence $\{\|\bar{\theta}_\ell-\tilde{x}\|^2\}_{\ell\geq0}$ converging a.s. to 0, which implies it also converges a.s. to 0. On the other hand, it follows from \eqref{eq:result1} that $\{\bar{\theta}_k\}_{k\geq0}$ is uniformly integrable. Then, by the dominated convergence theorem, we have $\lim_{k\to\infty}\mathbf{E}[\|\bar{\theta}_k-\tilde{x}\|] = 0$. Therefore, we conclude $\{\bar{\theta}_k-\tilde{x}\}_{k\geq0}$ converges to 0 with probability 1 and in mean. Noting that $\{x^i_k-\bar{\theta}_k\}_{k\geq0}$ converges to 0 with probability 1 and in mean from Theorem~\ref{theorem:consensus}-(2), we then have $\{x^i_k\}_{k\geq0}$ converges to $\tilde{x}$ with probability 1 and in mean, which completes the proof of part (1).

For part (2), due to the continuity of $f_\mu$, it holds that $\{f_\mu(x^i_k)\}_{k\geq0}$ converges a.s. to $f_\mu^\star$. To show $\lim_{k\to\infty}\mathbf{E}[f_\mu(x^i_k)]=f_\mu^\star$ by the dominated convergence theorem, it suffices to prove $\{f_\mu(x^i_k)-f_\mu^\star\}_{k\geq0}$ is uniformly integrable. One sufficient condition to ensure uniform integrability is that $\sup_{k\geq0}\mathbf{E}[|f_\mu(x^i_k)-f_\mu^\star|^2]<\infty$, see Sec.~4.12\cite{Siegrist2020}. Since 
\begin{align*}
\sup_{k\geq0}\mathbf{E}[|f_\mu(x^i_k)-f_\mu(\tilde{x})|^2]&\leq2\sup_{k\geq0}\mathbf{E}[|f_\mu(x^i_k)-f_\mu(\bar{\theta}_k)|^2]+2\sup_{k\geq0}\mathbf{E}[|f_\mu(\bar{\theta}_k)-f_\mu(\tilde{x})|^2]\\
&\leq2\hat{D}_\mu^2\sup_{k\geq0}\mathbf{E}[\|x^i_k-\bar{\theta}_k\|^2]+2\hat{D}_\mu^2\sup_{k\geq0}\mathbf{E}[\|\bar{\theta}_k-\tilde{x}\|^2]<\infty
\end{align*}
where the second inequality is due to the Lipschitz continuity of $f_\mu$ (c.f.~Lemma~\ref{lemma:property_f_mu}), and the last inequality follows from \eqref{eq:E_z-E_b_sq} and \eqref{eq:result1}, we thus conclude $\lim_{k\to\infty}\mathbf{E}[f_\mu(x^i_k)]=f_\mu^\star$. Likewise, we can also show that $\lim_{k\to\infty}\mathbf{E}[f(x^i_k)]$ exists.
Finally, for any $x^\star\in\mathcal{X}^\star$, applying Lemma~\ref{lemma:property_f_mu}-(1), the desired result follows from that
\begin{align*}
\lim_{k\to\infty}\mathbf{E}[f(x^i_k)] &\leq \lim_{k\to\infty}\mathbf{E}[f_\mu(x^i_k)]+ n\sqrt{m}{\mu}\hat{D}= f^\star_\mu+ n\sqrt{m}{\mu}\hat{D}\\
&\leq f_\mu(x^\star)+ n\sqrt{m}{\mu}\hat{D} \leq f^\star + 2n\sqrt{m}{\mu}\hat{D}.
\end{align*}
The proof is completed.
\end{Proof}

\begin{Remark}
As can be seen from Theorem~\ref{theorem:optimality}, the state information of all agents $x^i_k$, $i\in\mathcal{V}$ will converge to an optimal solution to the smoothed problem $x^\star_\mu\in\mathcal{X}^\star_\mu$ with probability 1 and in mean, and its global cost value will converge to a small neighborhood of the optimal value in mean with an error bound depending on the smoothing parameter $\mu$. To futher measure the distance between each agent's iterate $x^i_k$ and an optimal solution to the original problem $x^\star\in\mathcal{X}^\star$, it boils down to quantify the gap between $x^\star_\mu$ and $x^\star$. Since both $x^\star_\mu$ and $x^\star$ may not be unique under Assumption~\ref{assumption_local_f_lipschitz}, deriving an upper bound for $\|x^\star_\mu-x^\star\|$ is generally difficult, and there is no pleasant way to make this argument in the existing literature. However, if the local and global cost functions of the original problem satisfy some stronger assumptions (\textit{e.g.}, differentiability, strong convexity, Lipschitz gradient), we are able to show that $\|x^\star_\mu-x^\star\|\leq O(\mu)$.
\end{Remark}

If the local and global cost functions satisfy the following assumption, we are able to quantify the gap between the optimal solutions to problems \eqref{eq:cost_function} and \eqref{eq:smooth_cost_function}, $\|x^\star_\mu-x^\star\|$.

\begin{Assumption}\label{assumption_local_f_lipschitz_smooth}
The global cost function $f$ is $\chi$-strongly convex. Each local cost function $f_i$ can be non-convex, and is assumed to be differentiable. Its gradient is Lipschitz continuous with a constant $L$,
\textit{i.e.,} $\forall x, y \in \mathbb{R}^m$, there exists a constant $L>0$ such that $\|\nabla f_i(x) - \nabla f_i(y)\|\leq L\|x-y\|$, $\forall i\in\mathcal{V}$.
\end{Assumption}

\begin{Lemma}\label{lemma:optimality_solution_gap}
Suppose Assumption~\ref{assumption_local_f_lipschitz_smooth} holds. Then, the optimal solutions to problems \eqref{eq:cost_function} and \eqref{eq:smooth_cost_function}, denoted by $x^\star$ and $x^\star_\mu$, are unique, and satisfy that $\|x^\star_\mu-x^\star\|\leq\frac{n(m+3)^{3/2}\rho L}{2(1-\sqrt{1-\rho\chi})}\mu$, where $\rho\in(0,\frac{\chi}{n^2L^2}]$ is a constant.
\end{Lemma}
\begin{Proof}
We first give some properties of the cost functions $f_\mu$, which are directly obtained from Sec.~2\cite{Nesterov2017}.

\begin{Lemma}\label{lemma:property_f_mu_smooth}
(see Sec.~2\cite{Nesterov2017}) Under Assumption~\ref{assumption_local_f_lipschitz_smooth}, we have that
\begin{enumerate}
\item function $f_\mu = \sum_{i=1}^n f_{i,\mu}$ is strongly convex with a constant $\chi_\mu\leq\chi$,
\item The gradient $\nabla f_{i,\mu}(x)$ is Lipschitz continuous with a constant $\hat{L} \leq L$, \textit{i.e.}, $\|\nabla f_{i,\mu}(x) - \nabla f_{i,\mu}(y)\| \leq \hat{L}\|x - y\|$. 
% Moreover, $\nabla f_{i,\mu}(x)$ pointwise converges to $\nabla f_i(x)$ when $\mu$ tends to 0, \textit{i.e.}, $\forall x\in\mathbb{R}^m$, $\lim_{\mu\to0^+}\nabla f_{i,\mu}(x) = \nabla f_i(x)$.
\item The gradient $\nabla f_{i,\mu}(x)$ holds that $\|\nabla f_{i,\mu}(x)- \nabla f_i(x)\| \leq \frac12(m+3)^{3/2}\mu L$.
\end{enumerate}
\end{Lemma}

It follows from Assumption~\ref{assumption_local_f_lipschitz_smooth} and Lemma~\ref{lemma:property_f_mu_smooth}-(1) that both $f(x)$ and $f_\mu(x)$ are strongly convex. Hence, the optimal solutions to problems \eqref{eq:cost_function} and \eqref{eq:smooth_cost_function}, denoted by $x^\star$ and $x^\star_\mu$, are unique.

Now, we define $F(x,\mu) = f_\mu(x): \mathbb{R}^m\times\mathbb{R}_{\geq0}\to\mathbb{R}$ to explicitly quantify the effect of parameter $\mu$ on $f_\mu(x)$. Then, the unique optimal solution to $f_\mu(x)$, $x^\star_\mu$, can be expressed as a parameter function of $\mu$, denoted by $x^\star(\mu)$. Since
\begin{align*}
F(x,0) = \frac1{\kappa}\int_{\mathbb{R}^m} \sum_{i=1}^nf_i(x)\exp\bigg({-\frac12\|\zeta\|^2}\bigg)d\zeta=f(x)\bigg[\frac1{\kappa}\int_{\mathbb{R}^m}\exp\bigg({-\frac12\|\zeta\|^2}\bigg)d\zeta\bigg]=f(x),
\end{align*}
we have that $x^\star(0)$ is also the unique optimal solution to $f(x)$, \textit{i.e.}, $x^\star(0)=x^\star$. Then, to characterize the bound on $\|x^\star_\mu-x^\star\|$, we aim to study $\|x^\star(\mu)-x^\star(0)\|$, which boils down to the Lipschitz property of function $x^\star(\mu)$ at $\mu=0$. Next, we will derive such property, following the ideas in Theorem~2.1\cite{Dafermos1988}.

Since $f(x)$ is differentiable, hence $\nabla f_\mu(x)$ is well-defined for $\mu\geq0$. Moreover, when $\mu = 0$, it follows from the definition of $f_\mu(x)$ that $\nabla f_0(x) = \frac1{\kappa}\int_{\mathbb{R}^m} \sum_{i=1}^n\nabla f_i(x)\exp({-\frac12\|\zeta\|^2})d\zeta = \nabla f(x)(\frac1{\kappa}\int_{\mathbb{R}^m}\exp({-\frac12\|\zeta\|^2})d\zeta) = \nabla f(x)$. Now, we define $\nabla F(x,\mu) = \nabla f_\mu(x): \mathbb{R}^m\times\mathbb{R}_{\geq0}\to\mathbb{R}^m$ to explicitly quantify the effect of parameter $\mu$ on $\nabla f_\mu(x)$. We show that $\nabla F(x,\mu)$ is (i) strongly monotone in $x$ for any $\mu\geq0$, (ii) Lipschitz continuous in $x$ for any $\mu\geq0$, and (iii) Lipschitz continuous in $\mu$ at $\mu=0$ for any $x\in\mathbb{R}^m$. 

(i) For $\mu>0$, we have $(x-y)^T(\nabla F(x,\mu)-\nabla F(y,\mu)) = (x-y)^T(\nabla f_\mu(x)-\nabla f_\mu(y)) \geq \chi\|x-y\|^2$, $\forall x,y\in\mathbb{R}^m$ due to strong convexity of $f_\mu$ (c.f. Lemma~\ref{lemma:property_f_mu_smooth}-(1)). For $\mu=0$, we have $(x-y)^T(\nabla F(x,0)-\nabla F(y,0)) = (x-y)^T(\nabla f(x)-\nabla f(y)) \geq \chi\|x-y\|^2$, $\forall x,y\in\mathbb{R}^m$ due to strong convexity of $f$ (c.f. Assumption~\ref{assumption_local_f_lipschitz_smooth}). Thus, we have $\nabla F(x,\mu)$ is $\chi$-strongly monotone in $x$ for any $\mu\geq0$, \textit{i.e.}, $\forall x,y\in\mathbb{R}^m$, we have
\begin{align}
(x-y)^T(\nabla F(x,\mu)-\nabla F(y,\mu)) \geq \chi\|x-y\|^2, \quad \forall \mu\geq0. \label{eq:strong_monotone}
\end{align}

(ii) For $\mu>0$, we have $\|\nabla F(x,\mu)-\nabla F(y,\mu)\| = \|\nabla f_\mu(x)-\nabla f_\mu(y)\|\leq nL\|x-y\|$, $\forall x,y\in\mathbb{R}^m$ due to Lipschitz continuity of $f_\mu$ (c.f. Lemma~\ref{lemma:property_f_mu_smooth}-(2)). For $\mu=0$, we have $\|\nabla F(x,0)-\nabla F(y,0)\| = \|\nabla f(x)-\nabla f(y)\|\leq nL\|x-y\|$, $\forall x,y\in\mathbb{R}^m$ due to Lipschitz continuity of $f$ (c.f. Assumption~\ref{assumption_local_f_lipschitz_smooth}). Thus, we have $\nabla F(x,\mu)$ is $\chi$-Lipschitz continuous in $x$ for any $\mu\geq0$, \textit{i.e.}, $\forall x,y\in\mathbb{R}^m$, we have
\begin{align}
\|\nabla F(x,\mu)-\nabla F(y,\mu)\| \leq nL\|x-y\|, \quad \forall \mu\geq0. \label{eq:lipschitz_cont_x}
\end{align}

(iii) To show that $\nabla F(x,\mu)$ is Lipschitz continuous in $\mu$ at $\mu=0$ for any $x\in\mathbb{R}^m$, that is to show $\|\nabla F(x,\mu) - \nabla F(x,0)\| \leq K|\mu-0| = K\mu$ for some constant $K>0$, which directly follows from Lemma~\ref{lemma:property_f_mu_smooth}-(3) with $K = \frac12n(m+3)^{3/2}L$. Hence, we have $\nabla F(x,\mu)$ is Lipschitz continuous in $\mu$ at $\mu=0$ for any $x\in\mathbb{R}^m$, \textit{i.e.}, $\forall x\in\mathbb{R}^m$, we have
\begin{align}
\|\nabla F(x,\mu) - \nabla F(x,0)\| \leq \frac12n(m+3)^{3/2}\mu L. \label{eq:lipschitz_cont_mu}
\end{align}

Now, consider the map $G(x,\mu) = x-\rho\nabla F(x,\mu)$, where $(x,\mu)\in\mathbb{R}^m\times\mathbb{R}_{\geq0}$. For $0<\rho \leq \frac{\chi}{n^2L^2}$,
\begin{align}
\|G(x,\mu) - G(y,\mu)\|^2&\leq\|x-y\|^2 - 2\rho(x-y)^T(\nabla F(x,\mu)-\nabla F(y,\mu))\nonumber\\
&\quad\quad+\rho^2\|\nabla F(x,\mu)-\nabla F(y,\mu)\|^2 \nonumber\\
&\leq(1-2\rho\chi+n^2\rho^2L^2)\|x-y\|^2\leq (1-\rho\chi)\|x-y\|^2, \label{eq:contraction_map}
\end{align}
where the second inequality is due to \eqref{eq:strong_monotone} and \eqref{eq:lipschitz_cont_x}. It follows from \eqref{eq:contraction_map} that the map $G(x,\mu)$ is a contraction with respect to $x$. By the Banach fixed point theorem, the map $G(x,\mu)$ has a unique fixed point $x(\mu)$. On the other hand, any fixed point $x(\mu)$ of the map $G(x,\mu)$ is a solution to $F(x,\mu)$ (since $x(\mu)=G(x(\mu),\mu)\iff\nabla F(x(\mu),\mu) = 0$). Thus, we deduce that the unique optimal solution to $F(x,\mu)$ is the unique fixed point of the map $G(x,\mu)$, \textit{i.e.}, $x^\star(\mu) = G(x^\star(\mu),\mu)$. Then, by \eqref{eq:contraction_map}, we have
\begin{align*}
\|x^\star(\mu) - x^\star(0)\| &= \|G(x^\star(\mu),\mu) - G(x^\star(0),0)\| \\
&\leq \|G(x^\star(\mu),\mu) - G(x^\star(0),\mu)\| + \|G(x^\star(0),\mu) - G(x^\star(0),0)\|\\
&\leq \sqrt{1-\rho\chi}\|x^\star(\mu) - x^\star(0)\| + \|G(x^\star(0),\mu) - G(x^\star(0),0)\|.
\end{align*}
Noting that the last term
\begin{align*}
\|G(x^\star(0),\mu) - G(x^\star(0),0)\| &= \|(x^\star(0) - \rho\nabla F(x^\star(0),\mu)) - (x^\star(0) - \rho\nabla F(x^\star(0),0))\|\\
&=\rho\|\nabla F(x^\star(0),\mu) - \nabla F(x^\star(0),0)\| \leq\frac12n(m+3)^{3/2}\mu\rho L,
\end{align*}
where the last inequality is due to \eqref{eq:lipschitz_cont_mu}. Combining the above two relations, we obtain the desired result.
\end{Proof}

\begin{Corollary}
Under Assumptions~\ref{assumption_graph}, \ref{assumption_local_f_lipschitz}, \ref{assumption_local_f_lipschitz_smooth}, and the conditions in Theorem~\ref{theorem:optimality}, for $i\in\mathcal{V}$, we have sequence $\{x^i_k\}_{k\geq0}$ converges a.s. to a small neighborhood of the unique solution $x^\star$, with the gap satisfying $\lim_{k\to\infty}\mathbf{E}[\|x^i_k - x^\star\|]\leq \frac{n(m+3)^{3/2}\rho L}{2(1-\sqrt{1-\rho\chi})}\mu$.
\end{Corollary}
\begin{Proof}
The conclusion directly follows from Theorem~\ref{theorem:optimality} and Lemma~\ref{lemma:optimality_solution_gap} by noting that $\lim_{k\to\infty}\mathbf{E}[\|x^i_k - x^\star\|]\leq \lim_{k\to\infty}\mathbf{E}[\|x^i_k - x^\star_\mu\|]+\|x^\star_\mu - x^\star\|=\|x^\star_\mu - x^\star\|\leq\frac{n(m+3)^{3/2}\rho L}{2(1-\sqrt{1-\rho\chi})}\mu$.
\end{Proof}

\section{Numerical Simulations}\label{sec:simulation}

We demonstrate the performance of our method with two numerical examples. We first consider a simple numerical example for the purpose of verifying the results we derived in the previous section. Then, we consider the case of privacy enhancement in distributed optimization, which has been mentioned in Sec.~\ref{sec:motivating_ex} \textit{Motivating Example}. 

\subsection{A Simple Example}
\begin{figure}
    \centering
    \subfloat[local cost function]{{\includegraphics[height=125pt]{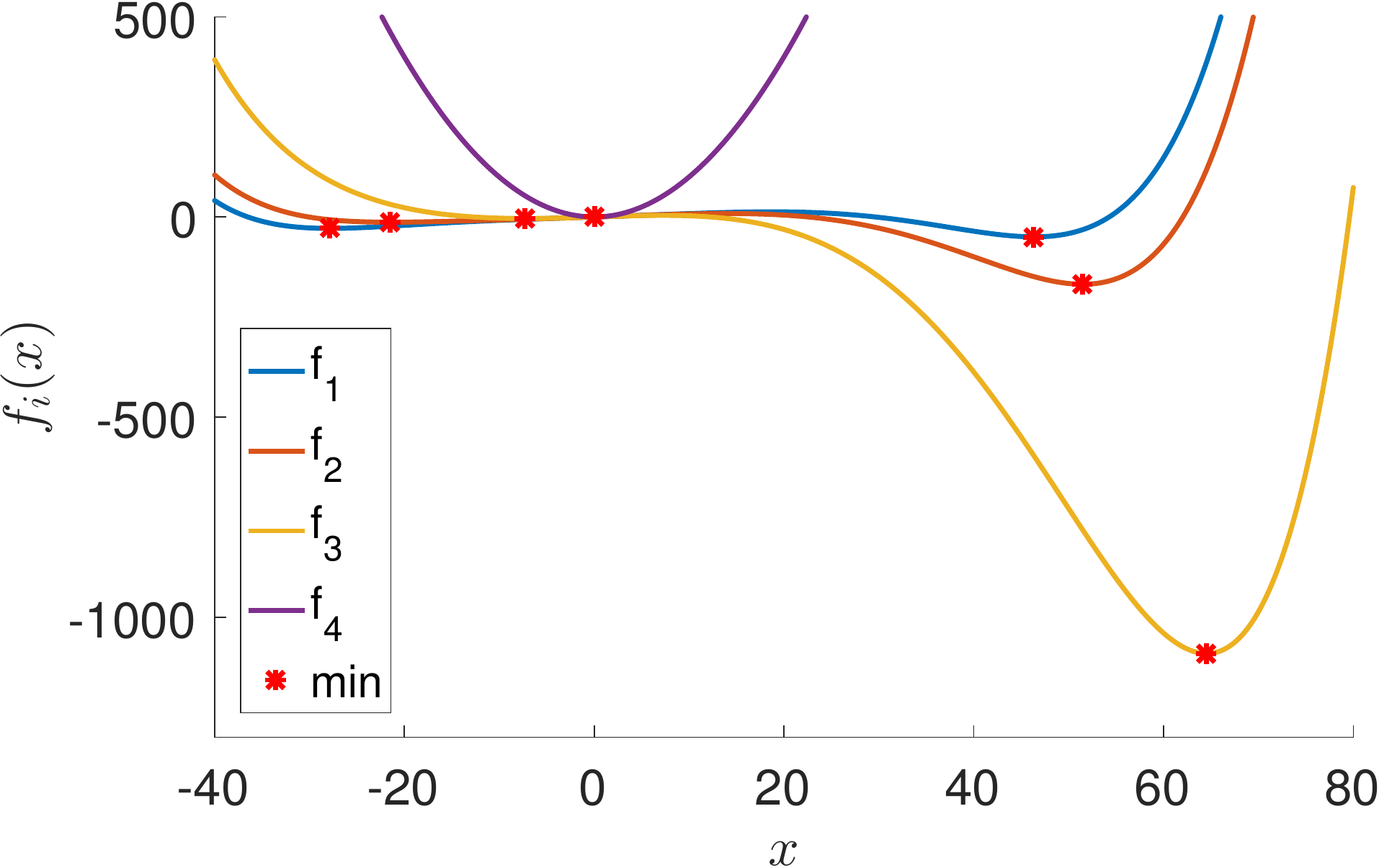} }}
    \subfloat[global cost function]{{\includegraphics[height=125pt]{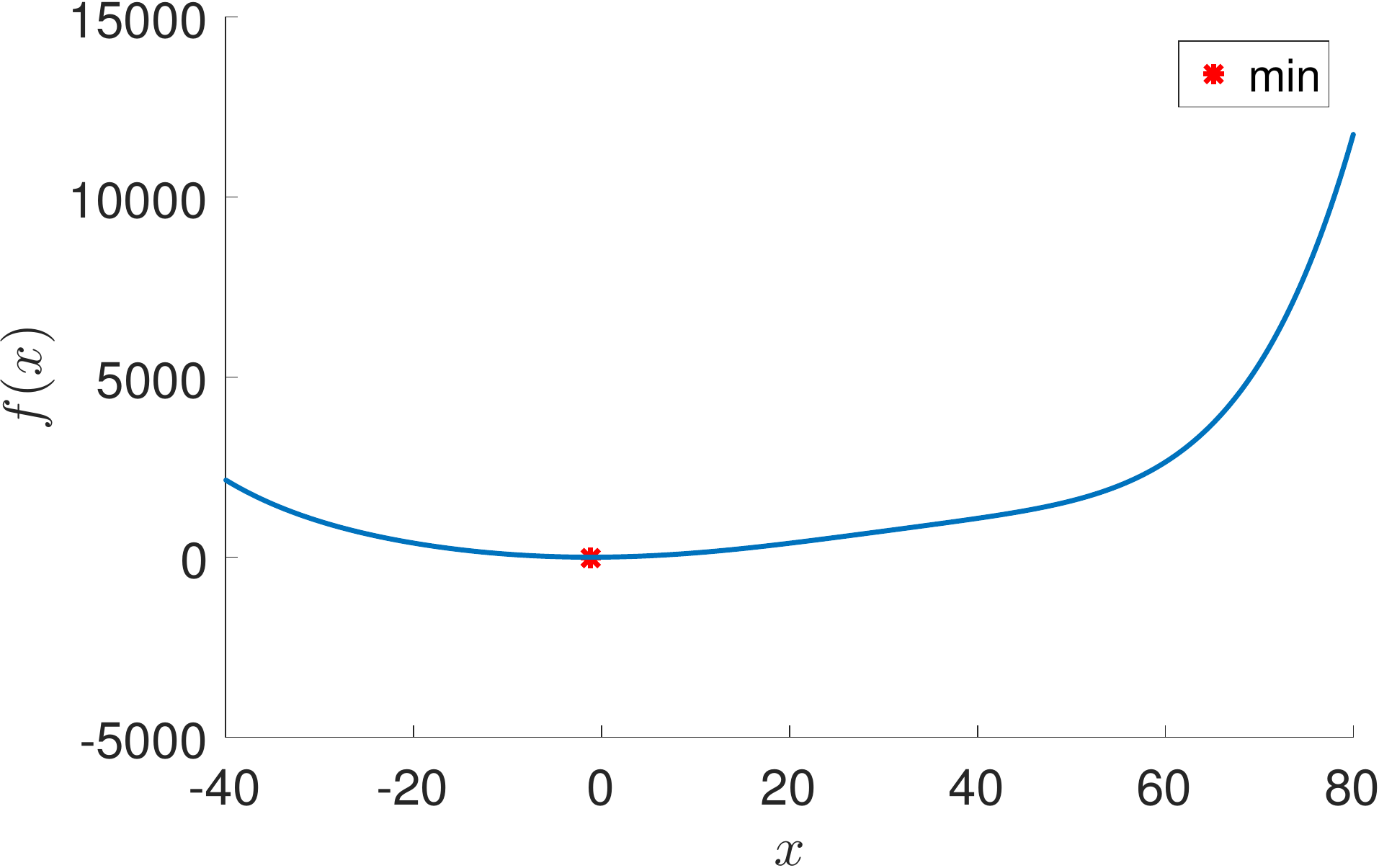} }}%
    \caption{Graphs of local cost function $f_i$ and global cost function $f$.}%
    \label{fig:local-global-cost.PNG}%
\end{figure}

\begin{figure}
\centering
\includegraphics[height=102pt]{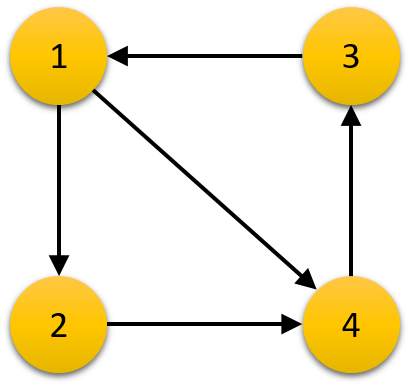}  % The printed column width is 8.4 cm.
\caption{Communication graph.}
\label{fig:network.png}
\end{figure}

Let the distributed optimization problem depicted in \eqref{eq:cost_function} be defined as follows
\begin{equation*}
\min f(x) = \sum_{i=1}^4 f_i(x),\quad x \in \mathbb{R},
\end{equation*}
where $f_1 = 3\times10^{-8}x^6 - 10^{-6}x^5 - 5\times10^{-5}x^4 + 5\times10^{-4}x^3 + 5\times10^{-3}x^2 + 0.8x$, $f_2 = 3\times10^{-8}x^6 + 10^{-6}x^5 - 5\times10^{-5}x^4 - 5\times10^{-4}x^3 + 5\times10^{-3}x^2 + 0.8x$, $f_3 = 3\times10^{-8}x^6 - 10^{-6}x^5 - 5\times10^{-5}x^4 - 5\times10^{-3}x^3 + 5\times10^{-3}x^2 + 0.8x$, and $f_4 = x^2$. The graphs of each individual cost $f_1, f_2, f_3, f_4$ and the global cost $f$ are plotted in Fig.~\ref{fig:local-global-cost.PNG}. Suppose the 4 agents are communicating through a network topology shown in Fig.~\ref{fig:network.png}.

Throughout the simulation, we let $[W_r]_{ij} = 1/|\mathcal{N}^{\text{in}}_i|$, $[W_c]_{ij} = 1/|\mathcal{N}^{\text{out}}_j|$, step-size $\alpha_k=1/(k+1)^{0.6}$, and $\epsilon = 0.1$. 

\begin{figure}
    \centering
    \subfloat[consensus ($x^i_k\to\bar{\theta}_k$)]{{\includegraphics[height=125pt]{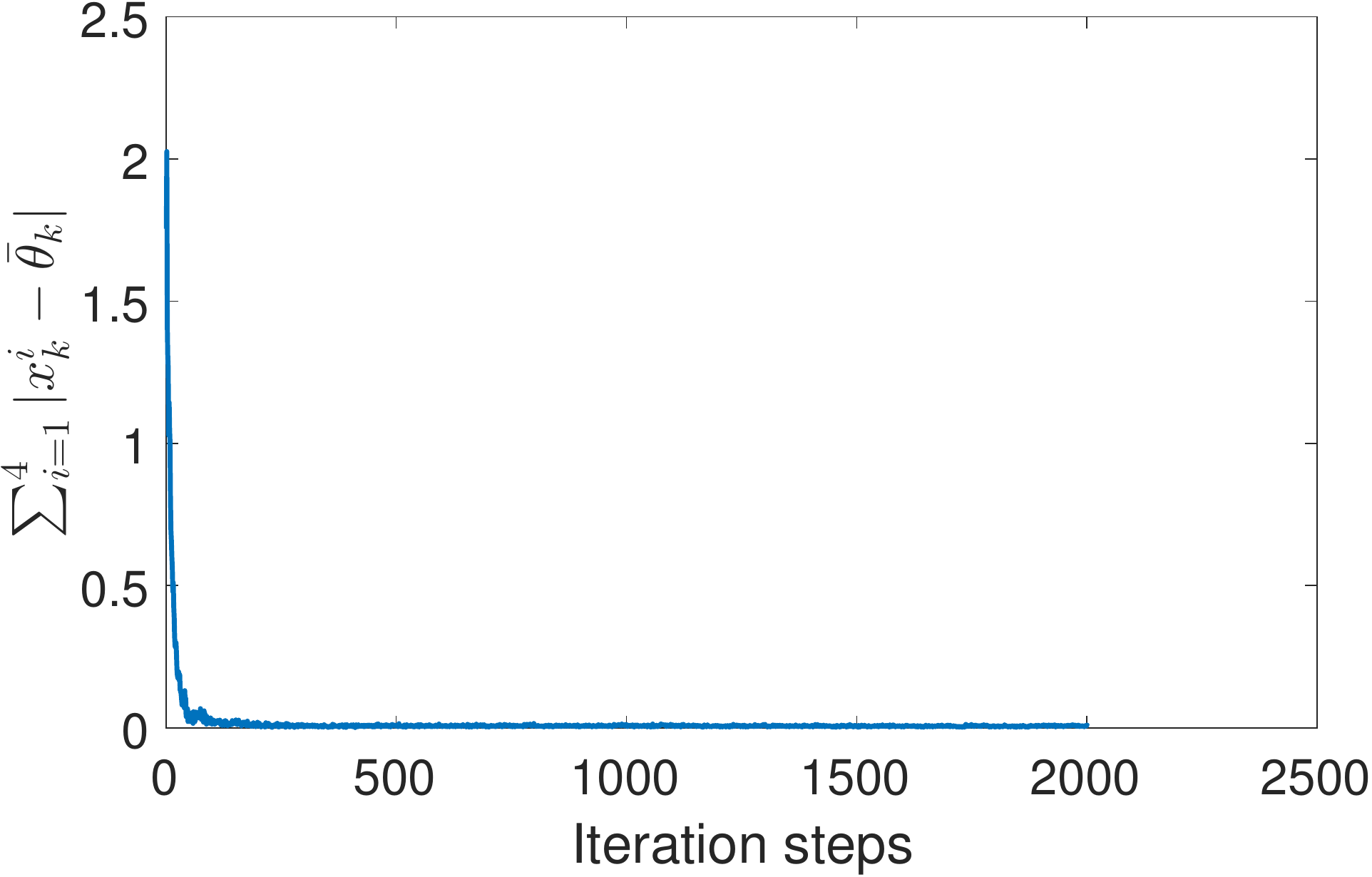} }}
    \subfloat[optimality ($f(\bar{\theta}_k)\to f^\star$)]{{\includegraphics[height=125pt]{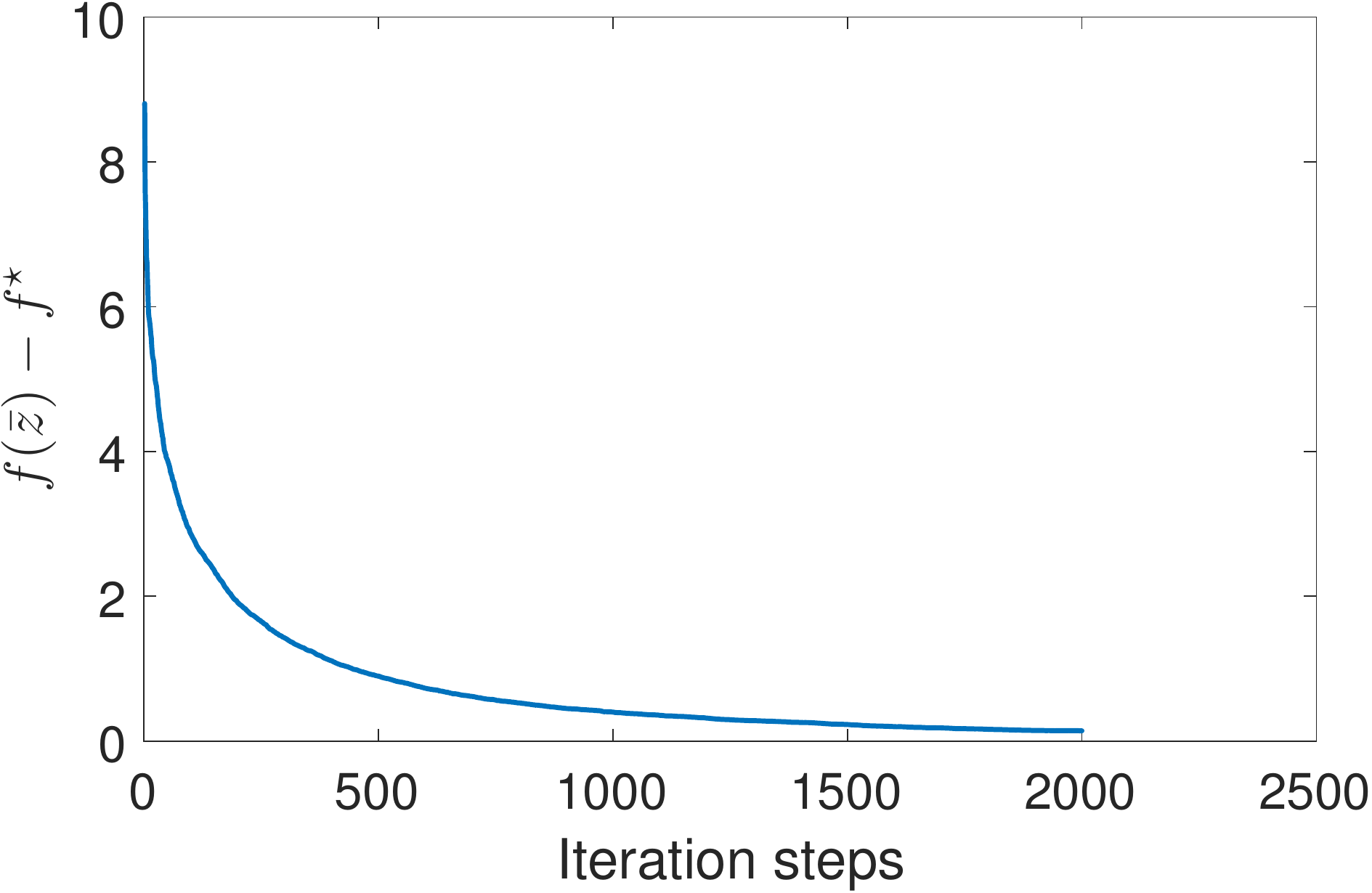} }}%
    \caption{Convergence performance of our method.}%
    \label{fig:convergence_of_x_and_f.PNG}%
\end{figure}

We first verify the convergence results of the proposed algorithm. Hence, we let $\mu = 10^{-2}$. The convergence results of both the consensus part and the optimality part, characterized by $\sum_{i=1}^4|x^i - \bar{\theta}|$ and $f(\bar{\theta}) - f^\star$, were plotted in Fig.~\ref{fig:convergence_of_x_and_f.PNG}-(a) and Fig.~\ref{fig:convergence_of_x_and_f.PNG}-(b), respectively. From both figures, it can be observed that $x^i_k$ converges to their average $\bar{\theta}_k$, whose function value $f(\bar{\theta}_k)$ also converges to the optimal value $f^\star$. These simulation results are consistent with the results derived in Theorems~\ref{theorem:consensus} and \ref{theorem:optimality}.

\begin{figure}
    \centering
    \subfloat[consensus]{{\includegraphics[height=125pt]{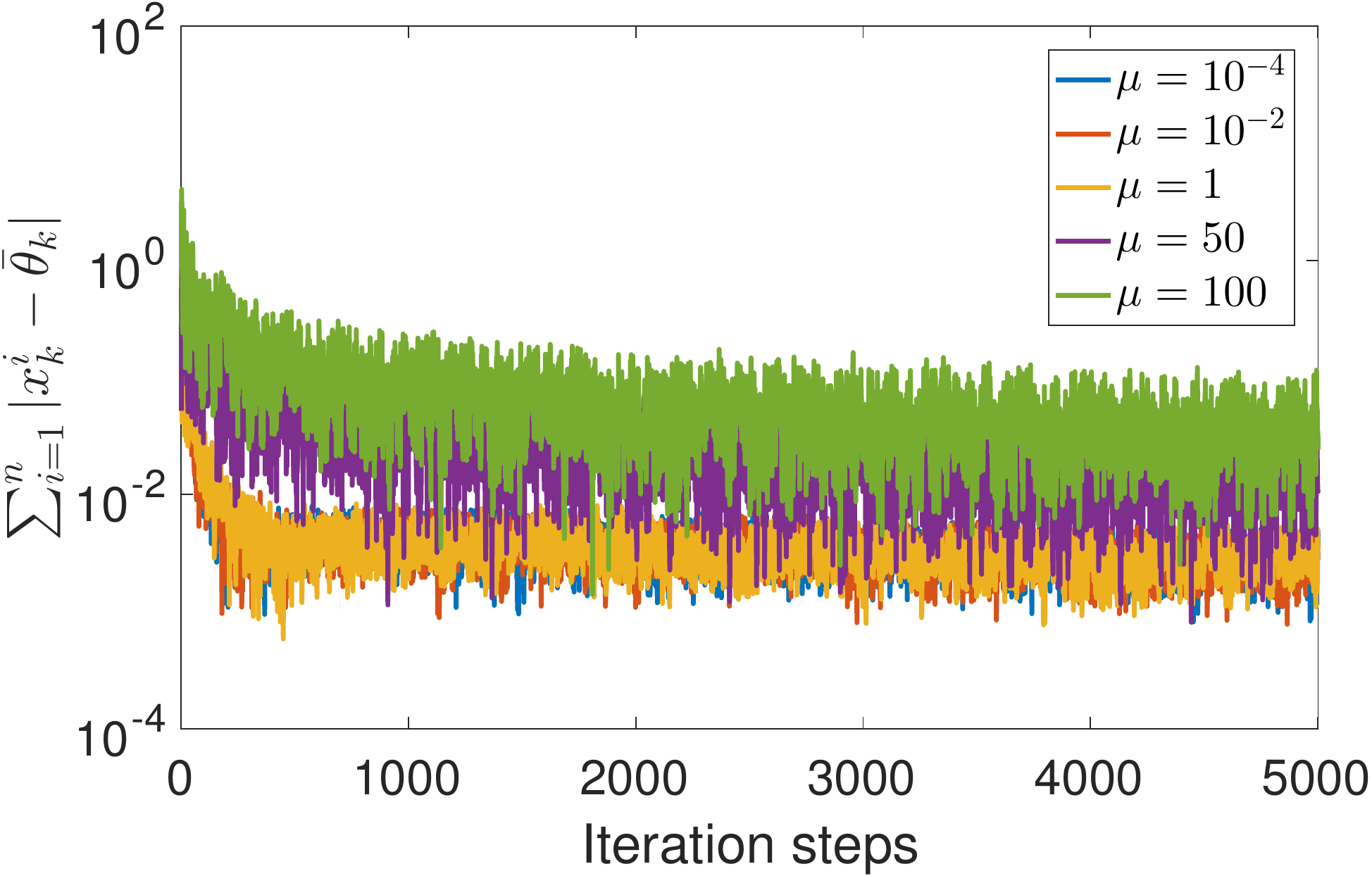} }}
    \subfloat[optimality]{{\includegraphics[height=125pt]{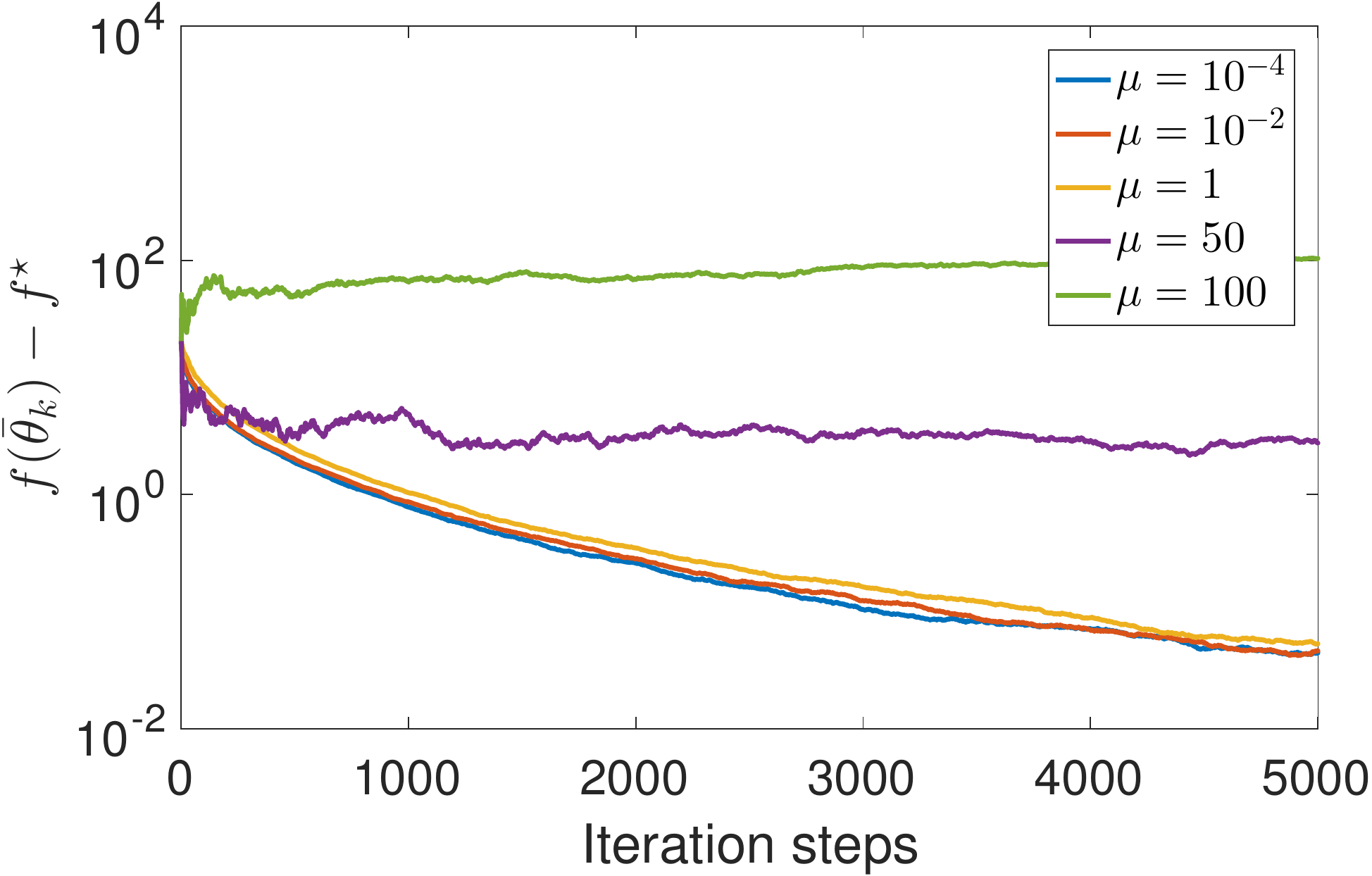} }}%
    \caption{Convergence results for different smoothing parameters $\mu$.}%
    \label{fig:smoothing_compare.PNG}%
\end{figure}

Then, we investigate the impact of the smoothing parameter $\mu$ on the convergence result of the proposed algorithm. Hence, in the simulation, $\mu$ was set to $10^{-4}$, $10^{-2}$, $1$, $50$ and $100$. The convergence results of both the consensus part and the optimality part under these five cases were plotted in Fig.~\ref{fig:smoothing_compare.PNG}-(a) and Fig.~\ref{fig:smoothing_compare.PNG}-(b), respectively. From both figures, it can be found that a smaller $\mu$ leads to both smaller consensus and optimality errors in general, which is due to the penalty term $2n\sqrt{m}{\mu}\hat{D}$. Overall, the influence on the convergence results is not significant with small $\mu$.

\subsection{Privacy Enhancement Example}
For the distributed optimization problem described in Sec.~\ref{sec:motivating_ex}
\begin{align*}
\min F(x) = \sum_{i=1}^n F_i(x),\quad x \in \mathbb{R}^m.
\end{align*}
We consider a logistic regression problem for binary classification over $n$ agents. Each agent $i$ is able to access $D_i$ data samples, $(a_{il},b_{il})\in\mathbb{R}^m\times\{-1,+1\}$, $l=1,\ldots,D_i$ where $a_{il}$ contains $m$-dimensional features, and $b_{il}$ is the corresponding binary label. The local cost function is given by $F_i(x) = \sum_{l=1}^{D_i}\ln[1+e^{-(x^Ta_{il})b_{il}}]+c_i\|x\|^2$, where $c_i$ is a regularization parameter.
For each agent $i$, functions $F_{i,j}, j \in \mathcal{N}^\text{out}_i$ are some fractional functions, \textit{i.e.}, $F_{i,j}(x) = \frac{d^i_j\|x\|^2}{1+\|x\|^2}$, whose coefficients $d^i_j$ are randomly generated. Then, each agent subtracts the self-generated functions for its out-neighbors $\sum_{j\in\mathcal{N}_i^{\text{out}}}F_{i,j}$, and combines the functions received from its in-neighbors $\sum_{s\in\mathcal{N}_i^{\text{in}}}F_{s,i}$ to form a new local cost function, given by
\begin{align*}
\tilde{F}_i(x) = F_i(x) - \sum_{j\in\mathcal{N}_i^{\text{out}}}F_{i,j} + \sum_{s\in\mathcal{N}_i^{\text{in}}}F_{s,i}.
\end{align*}
Hence, $\tilde{F}_i(x)$ may be non-convex due to the summation and subtraction of arbitrary functions. We apply the proposed algorithm to solve the new distributed optimization problem given by
\begin{align*}
\min \sum_{i=1}^n \tilde{F}_i(x),\quad x \in \mathbb{R}^m,
\end{align*}
which is equivalent to the original problem in terms of the optimal solution, but equipped with privacy enhancement settings.

Throughout the simulation, we let the smoothing parameter $\mu = 10^{-2}$ and the step-size $\alpha_k=1/(k+1)^{0.5}$. The convergence performance of the proposed algorithm characterized by the following two metrics: the consensus error given by $\sum_{i=1}^n\|x^i - \bar{\theta}\|$, and the optimality error defined as $F(\bar{\theta}) - F^\star$. 

\begin{figure}
    \centering
    \subfloat[consensus error]{{\includegraphics[height=125pt]{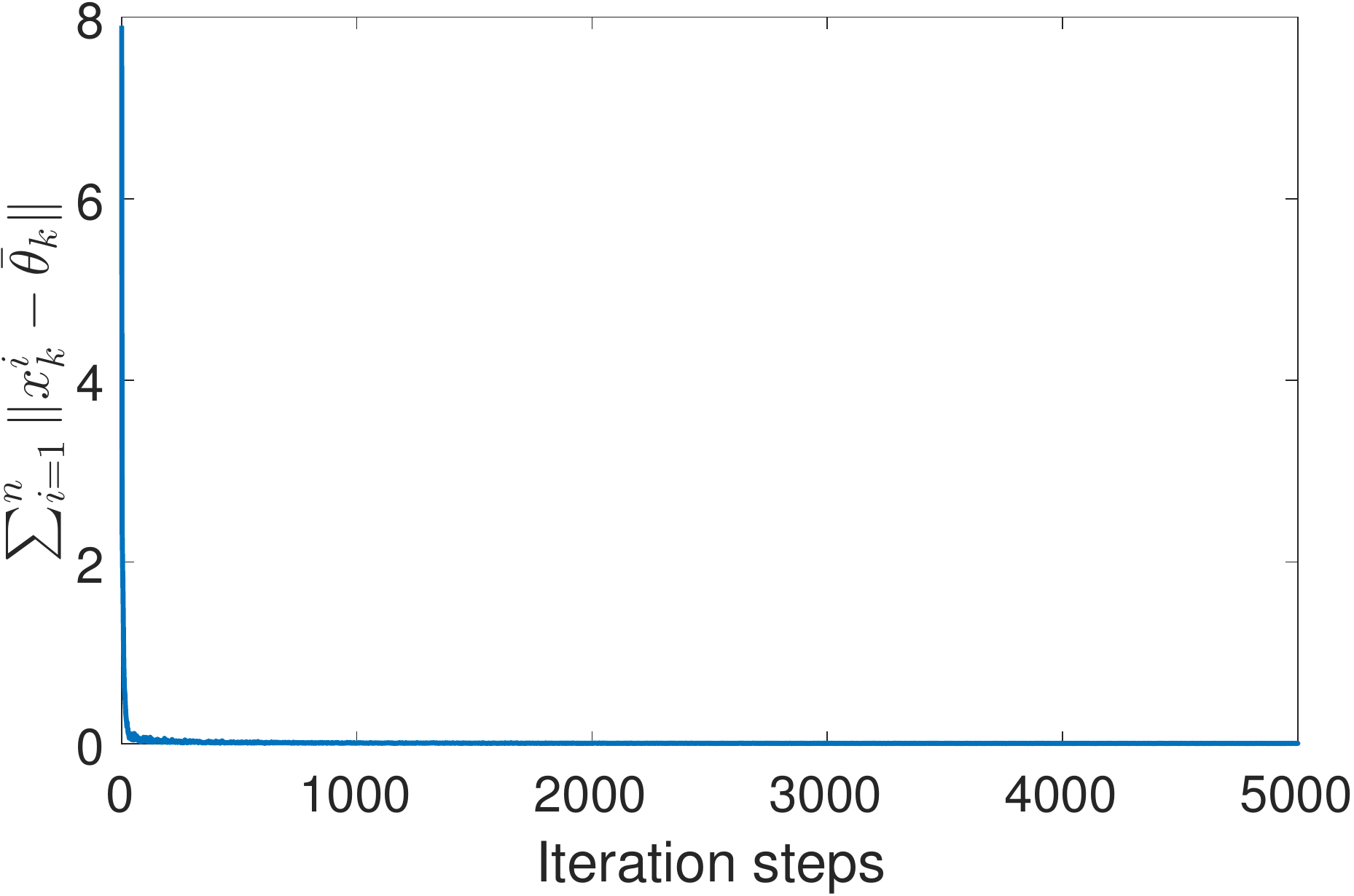} }}
    \subfloat[optimality error]{{\includegraphics[height=125pt]{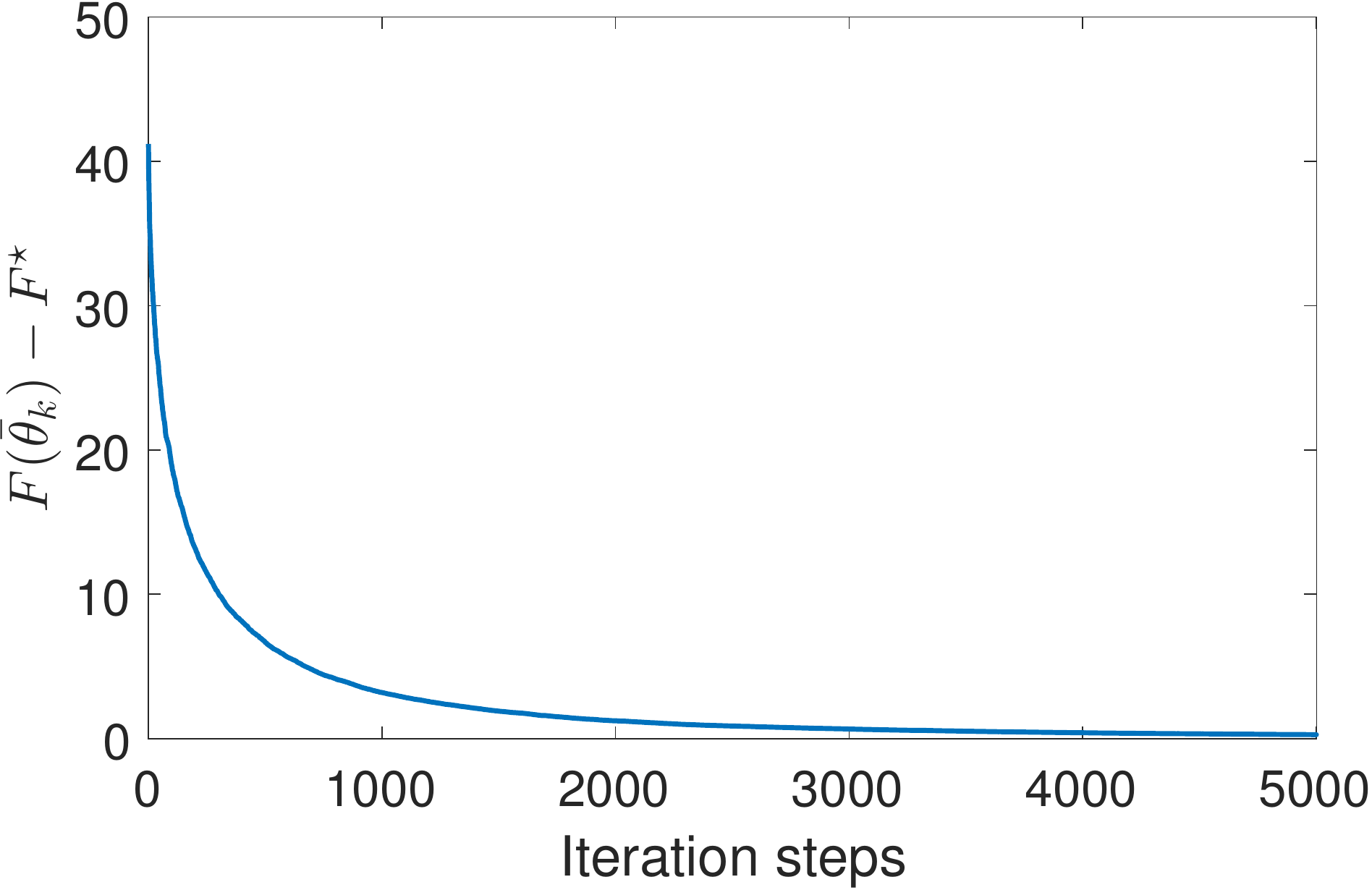} }}%
    \caption{Convergence performance of our method.}%
    \label{fig:convergence_of_x_and_f_2.PNG}%
\end{figure}

In the first simulation, we consider the dimension $m=2$, and the number of agents $n = 4$ with the communication graph as shown in Fig.~\ref{fig:network.png}. The simulation results were plotted in Fig.~\ref{fig:convergence_of_x_and_f_2.PNG}-(a) and Fig.~\ref{fig:convergence_of_x_and_f_2.PNG}-(b). From both figures, it can be observed that the proposed algorithm is able to achieve the convergence to the optimal solution.

\begin{figure}
    \centering
    \subfloat[consensus error]{{\includegraphics[height=125pt]{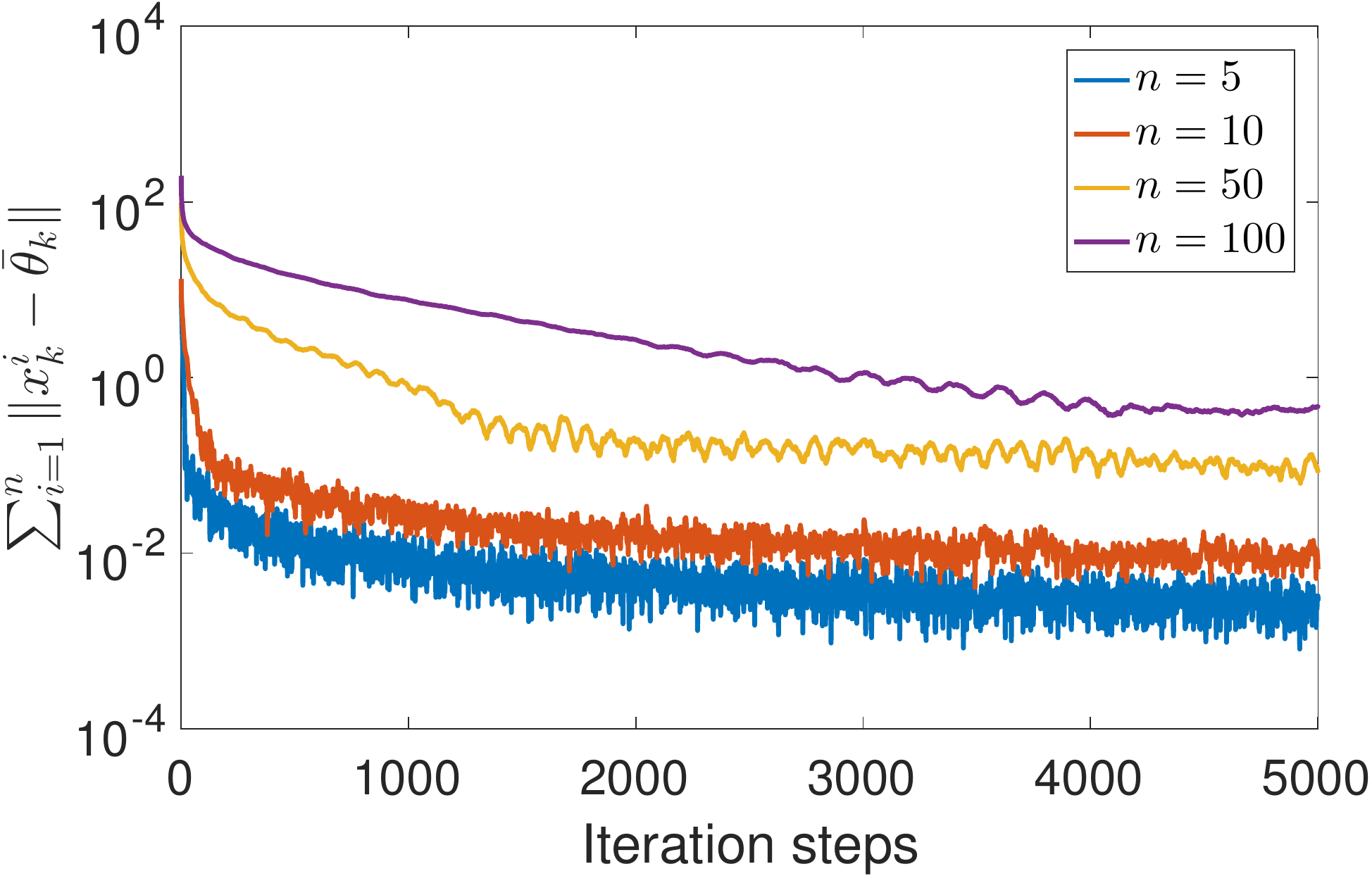} }}
    \subfloat[optimality error]{{\includegraphics[height=125pt]{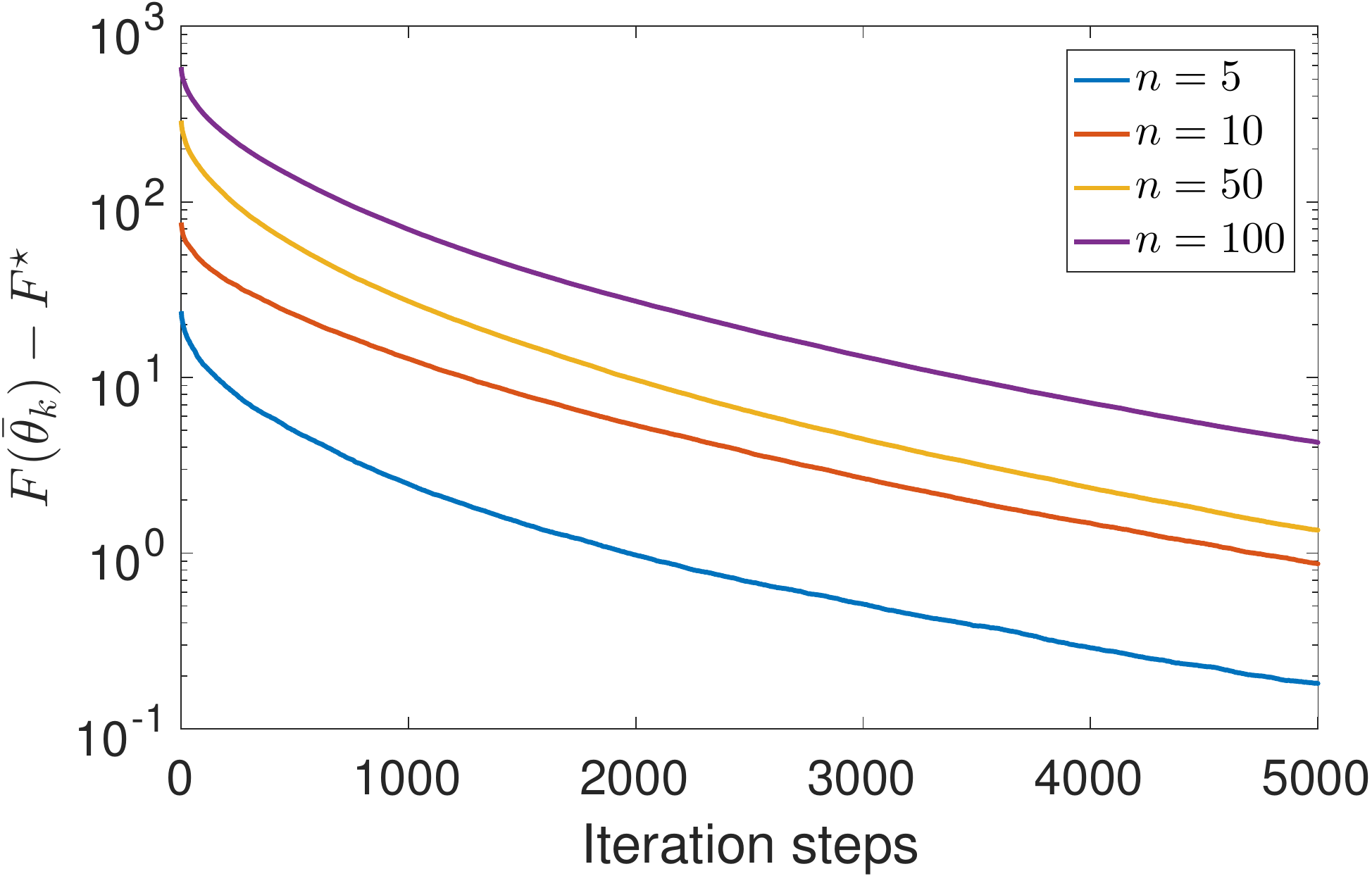} }}%
    \caption{Convergence results for different number of agents.}%
    \label{fig:numAgents_compare.PNG}%
\end{figure}

Next, we study the impact of the number of agents on the convergence results of our method. It is noted that the agents in large-scale applications may refer to sensors, computing units, decision makers, \textit{etc}. The number of agents is likely to vary from case to case. Hence, we set the number of agents $n = 5, 10, 50$ and $100$. The communication graph is supposed to be a circle with $\epsilon=10^{-6}$. Then, each agent have one in-neighbor and one out-neighbor, and hence the corresponding weight is set to $0.5$. The convergence results of both the consensus part and the optimality part under these four cases were plotted in Fig.~\ref{fig:numAgents_compare.PNG}-(a) and Fig.~\ref{fig:numAgents_compare.PNG}-(b), respectively. From both figures, it can be noticed that larger number of agents results in more discrepancies between agents' decisions during the iteration, hence leading to a longer time to converge to the optimal solution.

\begin{figure}
    \centering
    \subfloat[consensus error]{{\includegraphics[height=125pt]{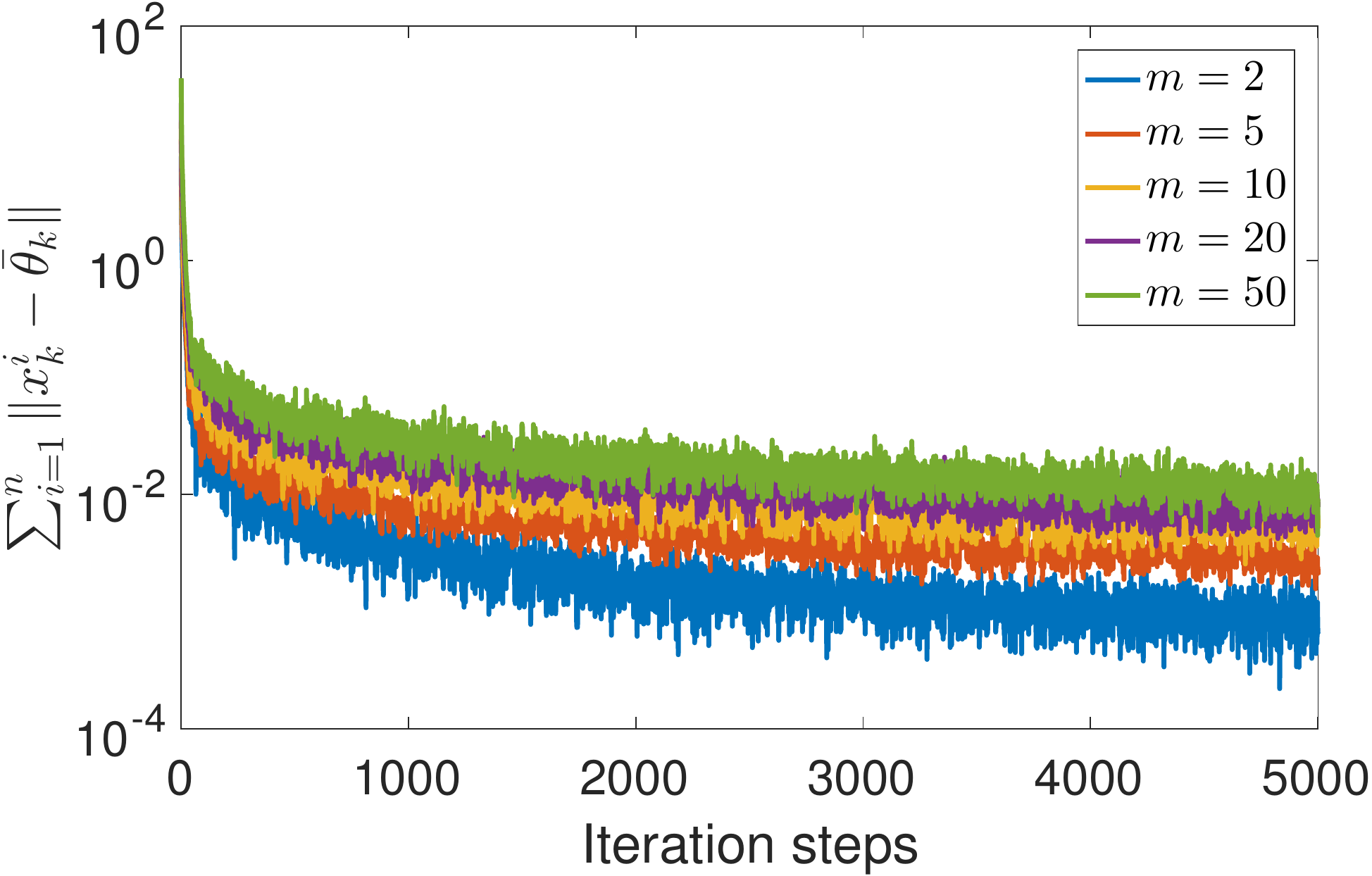} }}
    \subfloat[optimality error]{{\includegraphics[height=125pt]{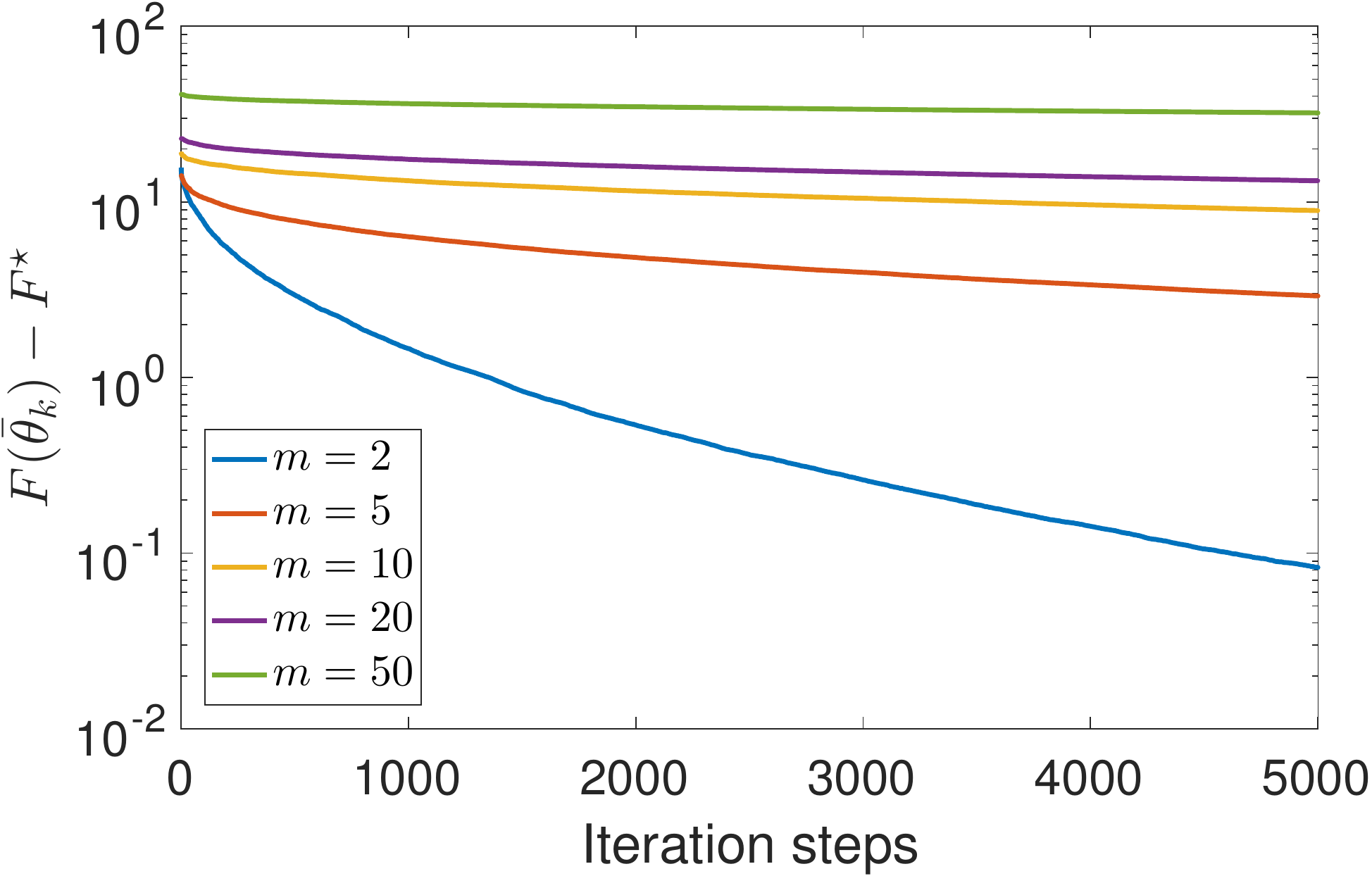} }}%
    \caption{Convergence results for different number of agents.}%
    \label{fig:nDim_compare.PNG}%
\end{figure}

Finally, we investigate the impact of the problem dimension on the convergence results of our method. Hence, we let $m=2$, $5$, $10$, $20$ and $50$. In this experiment, we consider the number of agents $n = 4$ with the communication graph as shown in Fig.~\ref{fig:network.png}. The simulation results were plotted in Fig.~\ref{fig:nDim_compare.PNG}-(a) and Fig.~\ref{fig:nDim_compare.PNG}-(b). As can be seen from both figures, both consensus and optimality errors become larger with the increase of the problem dimension.

\section{Conclusions}\label{sec:conclusion}
In conclusion, a randomized gradient-free distributed algorithm has been developed to solve a special type of distributed optimization problems where the global cost function is convex but each individual local cost function can be non-convex. 
We have proved that each agent's iterate approximately converges to the optimal solution both with probability 1 and in mean, and provided an upper bound on the optimality gap, characterized by the difference between the functional value of the iterate and the optimal value.
The convergence property of this algorithm has been analyzed in details. Finally, we have justified the performance of our method through a numerical example and an application in privacy enhancement. 
% Future work may include the analysis of the expected rate of convergence, the extension to a constrained optimization problem, the asynchronous update version, the time-varying communication network, \textit{etc.}

% \section*{Acknowledgments}
% This research was supported by Singapore Ministry of Education Academic Research Fund Tier 1 RG180/17(2017-T1-002-158).

% \subsection*{Author contributions}

% This is an author contribution text. This is an author contribution text. This is an author contribution text. This is an author contribution text. This is an author contribution text. 

% \subsection*{Financial disclosure}

% None reported.

% \subsection*{Conflict of interest}

% The authors declare no potential conflict of interests.

% \section*{Supporting information}

% The following supporting information is available as part of the online article:

% \noindent
% \textbf{Figure S1.}
% {500{\uns}hPa geopotential anomalies for GC2C calculated against the ERA Interim reanalysis. The period is 1989--2008.}

% \noindent
% \textbf{Figure S2.}
% {The SST anomalies for GC2C calculated against the observations (OIsst).}

% \appendix

% \section{Section title of first appendix\label{app1}}

% \nocite{*}% Show all bib entries - both cited and uncited; comment this line to view only cited bib entries;
\bibliographystyle{IEEEtran}
\bibliography{rgf_d_dgd_ncvx_reference}%

% \clearpage

% \section*{Author Biography}

% \begin{biography}{\includegraphics[width=66pt,height=86pt,draft]{empty}}{\textbf{Author Name.} This is sample author biography text this is sample author biography text this is sample author biography text this is sample author biography text this is sample author biography text this is sample author biography text this is sample author biography text this is sample author biography text this is sample author biography text this is sample author biography text this is sample author biography text this is sample author biography text this is sample author biography text this is sample author biography text this is sample author biography text this is sample author biography text this is sample author biography text this is sample author biography text this is sample author biography text this is sample author biography text this is sample author biography text.}
% \end{biography}

\end{document}